\documentclass[a4paper,11pt]{amsart}
\usepackage{amsmath}
\usepackage{amsfonts}
\usepackage{amssymb}
\usepackage{xcolor}
\usepackage[cp1250]{inputenc}
\usepackage[T1]{fontenc}

\newtheorem{theorem}{Theorem}[section]

\newtheorem{lemma}[theorem]{Lemma}
\newtheorem{corollary}[theorem]{Corollary}
\newtheorem{remark}[theorem]{Remark}
\newtheorem{conjecture}[theorem]{Conjecture}

\newcommand{\dow}[1][]{\textbf{Proof{#1}. }}
\newcommand{\tit}[1]{\textit{#1}}
\newcommand{\rightBox}{\hfill\Box}
\newcommand{\notextbf}{\normalfont}

\DeclareMathOperator{\mdeg}{mdeg}

\DeclareMathOperator{\sqrf}{sqrf}

\def\nl{[}
\def\np{]}

\makeatletter
\def\@cite#1#2{\nl{#1\if@tempswa , #2\fi}\np}
\makeatother

\makeatletter
\@namedef{subjclassname@2020}{%
	\textup{2020} Mathematics Subject Classification}
\makeatother

\begin{document}
\author{Daria Holik, Marek Karaś}
\thanks{\textit{Corresponding author:} Daria Holik, \textit{e-mail}: holikd@agh.edu.pl}
\title{Dependence of Homogeneous Components of Polynomials with Small Degree of Poisson Bracket}

\keywords{Poisson bracket, degree of Poisson bracket, polynomial ring, Jacobian determinant} 

\subjclass[2020]{Primary 13F20, 14R10; Secondary 16W20}


\maketitle

\pagestyle{myheadings}

\markboth{Daria Holik, Marek Karaś}{Dependence of Homogeneous Components}


\begin{abstract}
	Let $F,G\in\mathbb{C}[x_1,\ldots,x_n]$ be two polynomials in $n$ variables $x_1,\ldots,x_n$ over the complex numbers field $\mathbb{C}.$
	In this paper, we prove that if the degree of the Poisson bracket $[F,G]$ is small enough then there are strict constraints for homogeneous components of these polynomials (see Section \ref{Sect_formula_for_Gj}). 
	We also prove that there is a relationship between the homogeneous components of the  polynomial $F$ of degrees $\deg F -1$ and $\deg F-2$ (see Section \ref{sect_dependence_Fs_Fs_minus_1}) as well some results about divisibility of the homogeneous component of degree $\deg F- 1$ (see Section \ref{sect_divisibility_Fs}). 
	Moreover we propose, possibly an appropriate, reformulation of the conjecture of Yu regarding the estimation of the Poisson bracket degree of two polynomials (see Conjecture \ref{Conj_Yu_new}).
\end{abstract}

\section{introduction}

The research of this paper was inspired by the result of two papers.
One of them is the result of the paper \textit{Polynomial mappings with Jacobian determinant of bounded degree} due to  F. Leon Pritchard \cite{Pritchard} and the second one is the result from the paper \textit{Degree estimate for commutators} due to Vesselin Drensky and Jie-Tai Yu \cite{Drensky_Yu} (see also \cite{Makar_Yu} and \cite{Kuroda3}).

In his paper, Pritchard considered the following problem. 
For a given polynomial mapping $(F,G): \mathbb{C}^2 \rightarrow \mathbb{C}^2$ with small degree of the Jacobian determinant, i.e. with small degree of the polynomial $\frac{\partial F}{\partial x}\frac{\partial G}{\partial y} - \frac{\partial F}{\partial y}\frac{\partial G}{\partial x},$ where $F,G\in \mathbb{C}[x,y]$ are there some relations between homogeneous components of $F$ and $G$.
The main result of the paper gives such a relation.

On the other hand, Drensky and Yu, in their article consider the lower bound for the degree of the Poisson bracket (in the commutative case) or of the degree of the commutator (in the non-commutative case) of two polynomials $F$ and $G$ in $n$ variables $x_1, \dots ,x_n.$
Their considered the conjecture which was formulated by Yu \cite{Yu}. 
This conjecture says that, with some additional assumptions, the following inequality holds:
\begin{equation*}
\deg[F,G] > \min\{\deg F, \deg G\}.
\end{equation*}
Unfortunately the conjecture was showed to be false and, in the paper \cite{Drensky_Yu}, there was presented the counterexample, to the conjecture of Yu, due to Makar-Limanov (see \cite[Example 1.1]{Drensky_Yu}). 
Fortunately this example has some bad property. 
For instance this property can not be possessed by two polynomials that are components of the same polynomial authomorphism of $\mathbb{C}^n$. 
To be precise, let $F=F_1+ \dots + F_d$ and $G=G_1+ \dots + G_N$, where $F_i$ and $G_j$ are homogeneous polynomials of degree $i$ and $j,$ respectively.
In the example of Makar-Limanov the linear components of both polynomials $F$ and $G$ are linearly dependend. 
In the case when $F$ and $G$ are components of the same polynomial authomorphism of  $\mathbb{C}^n$ the linear components of $F$ and $G$ must be lineary independent. 
Thus, it is natural to state the following, new version of the conjecture of Yu:

\begin{conjecture}\label{Conj_Yu_new}
Let $F$ and $G$ be algebraically independent polynomials in $\mathbb{C}[x_1,\ldots,x_n].$ 
Let $F$ and $G$ generate their own centralizers in  $\mathbb{C}[x_1,\ldots,x_n],$ respectively. 
Suppose, also,  that $\deg F\nmid \deg G,$ $\deg G\nmid \deg F,$ $F(0)=G(0)=0$ and that linear parts of $F$ and $G$ are linearly independent. 
Then $\deg[F,G] > \min\{\deg F, \deg G\}$ or maybe some other (weaker) inequality holds, for example, $\deg[F,G] > \ c\cdot\min\{\deg F, \deg G\}$ for some constant $c$ which is not dependent of $F$ and $G.$ 
\end{conjecture}

Let us notice that the assumption that the homogeneous components of maximal degree of $F$ and $G$ are algebraically dependent, which appears in the original conjecture of Yu, is not necessary, because if the homogeneous components of maximal degree of $F$ and $G$ are algebraically independent then $\deg [F,G]=\deg F+\deg G>\min\{ \deg F,\deg G \}.$

The problem of estimation from below of the degree of the Poisson bracket $[F,G]$ is, of course, connected with the results about estimation of the lower bound of the degree of nonconstant elements of subalgebras generated by two polynomials.
The first such an estimation was given by Shestakov and Umirbaev \cite{shUmb}.
Their used this estimation and the theorem saying that each tame automorphism of $\mathbb{C}^3$ admits an elementary reduction or a reduction of some other types, to show that the famous Nagata example
\begin{equation}
\sigma:\mathbb{C}^3\ni (x,y,z)\mapsto (x+2y(y^{2}+zx)-z(y^{2}+zx)^{2},y-z(y^{2}+zx),z)\in \mathbb{C}^3
\end{equation}
is wild.
For information about four types of reductions defined by Shestakov and Umirbaev see for example \cite{shUmb1}.
In the paper \cite{Kuroda2} Shigeru Kuroda showed that there is no tame automorphism that admits a reduction of type IV.

Later in 2009 Leonid Makar-Limanov and Jie-Tai Yu  \cite{Makar_Yu} make an improvement of the estimation due to Shestakov and Umirbaev.
The improved estimation is the following
\begin{equation}
\deg P(f,g)\geq D(f,g)w_{\deg f,\deg g}P,
\end{equation} 
where $w_{\deg f,\deg g}P$ for polynomial $P=P(x,y)$ in two variables $x$ and $y$ denotes the weighted degree of $P$ such that $w_{\deg f,\deg g}x=\deg f$ and $w_{\deg f,\deg g}y=\deg g,$
and 
\begin{equation}
D(f,g)=1-\frac{\gcd(\deg f,\deg g)-(\deg(fg)-\deg [f,g])}{\deg f\deg g}.
\end{equation}
In the paper \cite{Makar_Yu} the authors used, instead of $\deg [f,g],$ the equivalent, in this context, notion of the degree of the following differential form $df\wedge dg.$
On the other hand in \cite{Kuroda3} Kuroda generalized Shestakov-Umirbaev inequality to the case of more than two polynomials.
In his article the degree of the Poisson bracket was replaced, which was very natural, by the degree of the differential form which is obtained as the wedge product of the differentials of these polynomials.

The numerous possible results about existence or nonexistence of tame automorphisms $P=(P_1,P_2,P_3):\mathbb{C}^3\rightarrow\mathbb{C}^3$ with a given multidegree of $P$ (by multidegree of $P$ we mean $\mdeg P=(\deg P_1,\deg P_2,\deg P_3)$) would be possible to be achieved if the Conjecture~\ref{Conj_Yu_new} appears to be true.
For some results about the multidegrees see for example \cite{EdoKanKarKur,Karas1,Karas3,Karas4,KarasZygad,SunChen}.
 
One of a possible way of attacking the above conjecture or maybe some weaker conjecture (with stronger assumptions) is to prove that small degree of $\deg[F,G]$ implies the relations between $F=F_1+ \dots + F_d$ and $G=G_1+ \dots + G_N$ (here, by $F_i$ and $G_j$ we denote the homogeneous component of $F$ of degree $i$ and the homogeneous component of $G$ of degree $j,$ respectively)
and then show that such a relation can not be satisfied. 
On the other hand such a relation can be viewed as a~generalization of the result of Pritchard.
Indeed, in the case $F, G \in \mathbb{C}[x,y]$ it is known that 
$\deg[F,G]=2+ \deg\left(\frac{\partial F}{\partial x}\frac{\partial G}{\partial y} - \frac{\partial F}{\partial y}\frac{\partial G}{\partial x}\right).$ 
Since in the paper, we study the pair of polynomials $F,G$ with small degree of $[F,G]$ we will assume that $\deg [F,G]<\deg F+\deg G.$ 
This implies that $[F_d,G_N]=0,$ and so by Lemma 2 in \cite{Umirbaev Yu} there exist $a,b\in \mathbb{C},$ $k_{1},k_{2}\in \mathbb{N}$ and a homogeneous polynomial $h$ such that $F_d=ah^{k_1}$ and $G_N=bh^{k_2}.$ 
Of course, withuot loss of generality, we can assume that $h$ is not a power of any other polynomial of lower degree, and that $a=1,$ because $\mathbb{C}$ is algebraically closed. 
Thus, in the rest of the paper, we will assume that $F_d=h^{k_1}$ and $G_N=bh^{k_2}$ with $b\in\mathbb{C}^{\star}=\mathbb{C}\setminus\{0\}$ denoted by $a_N.$
Of course, in this setup, we have $k_1=\frac{d}{\deg h}$ and $k_2=\frac{N}{\deg h}.$

The paper is organized as follows.
In  Section~\ref{Sect_Poisson_Brackets}, we recall the notion of Poisson bracket of two polynomials, formulate some properties that can be used in the calculations and present the main tool of the paper the so-called $H$-reduction method. 
In Section~\ref{Sect_preparatory_lemma}, we give some auxilliary lemmas that can be used in the next section in computation.
Section~\ref{Sect_formula_for_Gj} is devoted for establishing the general formulas for the homogeneous components of the polynomial $G$ in terms of the homogeneous components of the polynomial $F.$ 
Then, in the following section we use the above mentioned general formulas to prove that if $\deg [F,G]$ is small enough then the homogeneous components of the polynomial $F$ of degree $\deg F-1$ is divisible by $h$ in the case of $\deg h\leq 2$ and by $\sqrf h$ in the case $\deg h>2$ (for definition of  $\sqrf h$ see Section \ref{sect_divisibility_Fs}), where $h$ is the polynomial such that the homogeneous components of maximal degree of $F$ and $G$ are proportional to some powers of $h.$
Under the assumption that $h$ is a square-free polynomial we obtain that if $\deg [F,G]$ is smaller and smaller then $F_{s}$ (the homogeneous component of $F$ of degree $s,$ where $s=\deg F-1$) is divisible by higher and higher power of $h.$
The next section give some other relations between homogeneous components of the polynomial $F$ in some special cases of $\deg F.$
More precisely, we show that there is a strong relation between $F_{s-1}$ and $F_s,$ but with some assumptions about the degree of $F$ and $G,$ and of course, about $\deg [F,G].$

Let us notice that the assumption that $h$ is square-free geometrically means that the components at infinity of $F=0$ (and $G=0$) all have the same multiplicity.
The notion of square-free polynomials seems to be naturally arrising in the context of Jacobian Conjecture or related topics (see for example \cite{Jedrzejewicz_Zielinski}).

At the end of the introduction let us notice that all given in the present article results are valid not only for $\mathbb{C}$ but for any algebraically closed field of characteristic zero.

\section{Poisson bracket and $H$-reduction method}\label{Sect_Poisson_Brackets}
In this section we develop the main tool that will be used in the next
sections of the paper.

 We start with the following lemma in which we use the
notation 
\begin{equation*}
\mathbb{C}[x_{1},\ldots ,x_{n}]_{d}=\left\{ f\in \mathbb{C}[x_{1},\ldots ,x_{n}]:f%
\text{ is homogeneous of degree }d\right\} \cup \{0\}.
\end{equation*}

\begin{lemma}
\label{Lemma_H_reduction}Let $H$ be a nonconstant homogeneous
polynomial that is not a power of any other polynomial of lower degree (in particular, if $H$ is a square-free polynomial)
 and let $P$ be any homogeneous polynomial such that 
\begin{equation*}
\left[ H,P\right] =0.
\end{equation*}
Then there exist $a\in \mathbb{C}$ and $k\in \mathbb{N}$ such that 
\begin{equation*}
P=aH^{k}.
\end{equation*}
Moreover, if $P\in \mathbb{C}[x_{1},\ldots ,x_{n}]_{d}$ and $\deg H\nmid d,$
then $P=0.$
\end{lemma}

For the convenience of the reader, before we give the proof of the above lemma, let us recall that for any $f,g\in \mathbb{C}[x_{1},\ldots ,x_{n}],$ by $\left[ f,g\right] $ we denote the Poisson bracket of $f$ and $g,$ which is the
following sum: 
\begin{equation*}
\sum_{1\leq i<j\leq n}\left( \frac{\partial f}{\partial x_{i}}\frac{\partial
g}{\partial x_{j}}-\frac{\partial f}{\partial x_{j}}\frac{\partial g}{%
\partial x_{i}}\right) \left[ x_{i},x_{j}\right] ,
\end{equation*}
where $[x_{i},x_{j}]$ can be viewed as formal objects satisfying the condition (see for example \cite{shUmb1})
\begin{equation*}
\lbrack x_{i},x_{j}]=-[x_{j},x_{i}]\qquad \text{for all }i,j.
\end{equation*}
We also define 
\begin{equation*}
\deg \left[ x_{i},x_{j}\right] =2\qquad \text{for all }i\neq j,
\end{equation*}
$\deg 0=-\infty $ and 
\begin{equation*}
\deg \left[ f,g\right] =\max_{1\leq i<j\leq n}\deg \left\{ \left( \frac{%
\partial f}{\partial x_{i}}\frac{\partial g}{\partial x_{j}}-\frac{\partial f%
}{\partial x_{j}}\frac{\partial g}{\partial x_{i}}\right) \left[
x_{i},x_{j}\right] \right\} .
\end{equation*}

Since $2-\infty =-\infty ,$ we have 
\begin{equation*}
\deg [f,g]=2+\underset{1\leq i<j\leq n}{\max }\deg \left( \frac{\partial f}{%
\partial x_{i}}\frac{\partial g}{\partial x_{j}}-\frac{\partial f}{\partial
x_{j}}\frac{\partial g}{\partial x_{i}}\right) .
\end{equation*}

From the above equality we have 
\begin{equation}
\deg \left[ f,g\right] \leq \deg f+\deg g
\label{deg_Poisson_bracket_deg_f_g}
\end{equation}
and $\deg [f,g]\geq 2$ iff $f$ and $g$ are algebraically independent.

\dow[ of Lemma \ref{Lemma_H_reduction}]
Since $\left[ H,P\right] =0,$ it follows that $H$ and $P$ are algebraically
dependent and so by Lemma 2 in \cite{Umirbaev Yu} there exist $a,b\in \mathbb{C}%
,$ $k_{1},k_{2}\in \mathbb{N}$ and a homogeneous polynomial $h$ such that 
\begin{equation*}
P=ah^{k_{1}}\qquad \text{and\qquad }H=bh^{k_{2}}.
\end{equation*}
Since $H$ is  not a power of any other polynomial of lower degree, we conclude that $k_{2}=1$ and so we can take $h=H.$
Thus, in particular, if $a\neq 0$ (i.e. $P\neq 0$), then $\deg P$ is
divisible by $\deg H.$ 

$\rightBox$

The mentioned $H$-reduction method is the statement given in the following corollary.

\begin{corollary}
\label{Corr_H_reduction}Let $H$ be a  nonconstant homogeneous
polynomial that is not a power of any other polynomial of lower degree and let $P$ be any polynomial such that 
\begin{equation*}
\left[ H,P\right] =0.
\end{equation*}
Then $P\in \mathbb{C}[H].$
\end{corollary}

\dow
Let $d=\deg P$ and let $P=P_{0}+\cdots +P_{d}$ be the homogeneous
decomposition of $P.$ Since $\left[ H,P\right] =0,$ it follows that 
\begin{equation}
\left[ H,P_{i}\right] =0\qquad \text{for }i=0,\ldots ,d.
\label{Row_H_reduction_1}
\end{equation}
In particular, $\left[ H,P_{d}\right] =0.$ Since $P_{d}\neq 0$ (by
definition of $d$), it follows that $d=k\deg H$ for some $k\in \mathbb{N}.$ By (%
\ref{Row_H_reduction_1}) and Lemma \ref{Lemma_H_reduction} there exist $%
a_{0},\ldots ,a_{k}\in \mathbb{C},$ $a_{k}\neq 0$ such that 
\begin{equation*}
P_{l\deg H}=a_{l}H^{l}\qquad \text{for }l=0,\ldots ,k
\end{equation*}
and $P_{i}=0$ for $i\notin \{0,\deg H,2\deg H,\ldots ,k\deg H\}.$
$\rightBox$

We will also use the following fact that is easy to check.

\begin{lemma}
Let $P\in \mathbb{C}[x_{1},\ldots ,x_{n}].$ For any $Q,R\in \mathbb{C}%
[x_{1},\ldots ,x_{n}]$ and $\alpha ,\beta \in \mathbb{C}$ we have 
\begin{eqnarray*}
\left[ P,QR\right] &=&Q\left[ P,R\right] +R\left[ P,Q\right] , \\
\left[ P,\alpha Q+\beta R\right] &=&\alpha \left[ P,Q\right] +\beta \left[
P,R\right] , \\
\left[ P,Q\right] &=&-\left[ Q,P\right] .
\end{eqnarray*}
In other words, the mappings $Q\mapsto \left[ P,Q\right] $ and $Q\mapsto
\left[ Q,P\right] $ are $\mathbb{C}$-derivations.
\end{lemma}

\section{Preparatory lemmas}\label{Sect_preparatory_lemma}
We assume that are given the following two polynomials of degree, respectively, $d$ and $N$ with $d\leq N:$
\begin{equation}
F=F_1+\cdots+F_d\qquad\text{and}\qquad G=G_1+ \cdots+G_{N},
\end{equation}
where $F_i$ and $G_j$ are homogeneous polynomials of degree $i$ and $j,$ respectively,
with 
\begin{equation}
F_d=h^{\frac{d}{t}}\qquad\text{and}\qquad G_N=a_Nh^{\frac{N}{t}},
\end{equation}
where $h$ is a homogeneous polynomial of degree $t$ and $a_N$ is a nonzero constant. Of course, we must have $t|d$ and $t|N.$


Let us notice that:
\begin{eqnarray}
\left[ F_i,F_1^{\alpha_1}\cdots F_s^{\alpha_s}h^{\alpha_0}\right] &=&
\alpha_0  F_1^{\alpha_1}\cdots F_s^{\alpha_s}h^{\alpha_0-1}[F_i,h]
+\sum_{j=1}^s\alpha_j F_1^{\alpha_1}\cdots F_s^{\alpha_s}F_j^{-1}h^{\alpha_0}[F_i,F_j]\label{row_A},
\end{eqnarray}
where $s=d-1.$

\begin{lemma}\label{Lem_sumGiGi1Gi2}
Assume that polynomials
\begin{equation}
G_u=\sum_{\substack{\alpha=(\alpha_0,\ldots,\alpha_s)\in\mathbb{Z}\times\mathbb{N}^s
}}
c^{(u)}_{\alpha}F_1^{\alpha_1}\cdots F_s^{\alpha_s}h^{\alpha_0}
\qquad u=i,\ldots,i+s-1,
\end{equation}
with $c^{(u)}_{\alpha}\in \mathbb{C},$
satisfy the following condition:\vspace{0.3cm}\newline
(C1) \hspace{1cm}$ (\alpha_k+1)c^{(i+s-l)}_{\alpha+e_k}=
(\alpha_l+1)c^{(i+s-k)}_{\alpha+e_l},$\hspace{0.5cm} for all $1\leq k,l\leq s$ and $\alpha\in\mathbb{Z}\times\mathbb{N}^s,$
\vspace{0.3cm}\newline
where $e_0=(1,0,\ldots,0),$ $e_1=(0,1,0,\ldots,0),$ $\ldots,$  $e_s=(0,\ldots,0,1).$ 
Then, the sum $[F_s,G_i]+[F_{s-1},G_{i+1}]+\cdots+[F_1,G_{i+s-1}]$ can be written as
\begin{eqnarray}\label{row_GiGi1Gi2}
&&\left[  \sum_{\alpha\in\mathbb{Z}\times\mathbb{N}_+\times\mathbb{N}^{s-1}}\tfrac{\alpha_0}{\alpha_1}
c^{(i+s-1)}_{\alpha-e_1}
F_1^{\alpha_1}\cdots F_s^{\alpha_s}h^{\alpha_0-1}
\right.  \\
&& + 
\sum_{\alpha\in\mathbb{Z}\times\{0\}\times\mathbb{N}_+\times\mathbb{N}^{s-2}}\tfrac{\alpha_0}{\alpha_2}
c^{(i+s-2)}_{\alpha-e_2}
F_2^{\alpha_2}\cdots F_s^{\alpha_s}h^{\alpha_0-1} \nonumber\\
&&
\vdots
\nonumber\\
&& \left. 
+
\sum_{\alpha\in\mathbb{Z}\times\{(0,\ldots,0)\}\times\mathbb{N}_+}\tfrac{\alpha_0}{\alpha_s}
c^{(i)}_{\alpha-e_s}
F_s^{\alpha_s}h^{\alpha_0-1}\quad , \quad h \right].\nonumber
\end{eqnarray}
\end{lemma}
\dow
By (\ref{row_A}) we see that $[F_s,G_i]+[F_{s-1},G_{i+1}]+\cdots+[F_1,G_{i+s-1}]$ is the sum of the following summands:
\begin{eqnarray}
\sum_{\alpha\in\mathbb{Z}\times\mathbb{N}^s}F_1^{\alpha_1}\cdots F_s^{\alpha_s}h^{\alpha_0}
\left\{ (\alpha_k+1)c^{(i+s-l)}_{\alpha+e_k}[F_l,F_k]+
(\alpha_l+1)c^{(i+s-k)}_{\alpha+e_l}[F_k,F_l]
\right\}\label{row_summFkFl}
\end{eqnarray}
for all $1\leq k<l\leq s$
and
\begin{eqnarray}
&&
\sum_{k=1}^s\sum_{\alpha\in\mathbb{Z}\times\mathbb{N}^s}\alpha_0 c^{(i+s-k)}_{\alpha} F_1^{\alpha_1}\cdots F_s^{\alpha_s}h^{\alpha_0-1}[F_k,h] \label{row_pochIloczynu}.
\end{eqnarray}

Let us notice that (\ref{row_summFkFl}) equals to zero by (C1).

We can assume that for all $u,$  $c^{(u)}_{(\alpha_0,\ldots,\alpha_s)}=0,$ if $\alpha_k<0$ for some $k=1,\ldots,s.$ Then, by renumbering of the summs in (\ref{row_pochIloczynu}), the summand (\ref{row_pochIloczynu}) can be rewritten as follows
\begin{eqnarray}
&\sum\limits_{\alpha\in\mathbb{Z}\times\mathbb{N}^s}\alpha_0 h^{\alpha_0-1} &
\left\{ c^{(i+s-1)}_{\alpha-e_1} F_1^{\alpha_1-1}F_2^{\alpha_2}\cdots F_s^{\alpha_s}[F_1,h] 
+ c^{(i+s-2)}_{\alpha-e_2} F_1^{\alpha_1}F_2^{\alpha_2-1}F_3^{\alpha_3}\cdots F_s^{\alpha_s}[F_2,h] 
\right. \\
&&
\left.
+\cdots+
c^{(i)}_{\alpha-e_s} F_1^{\alpha_1}\cdots F_{s-1}^{\alpha_{s-1}}F_s^{\alpha_s-1}[F_s,h]
\right\} .\nonumber
\end{eqnarray} 
Assume that $\alpha_k\neq 0$ for some $k\in\{ 1,\ldots,s\}$ and set $e_0=(1,0,\ldots,0),$ $e_1=(0,1,0,\ldots,0),$ $\ldots,$  $e_s=(0,\ldots,1).$ 
Then, one can see that
\begin{eqnarray}
&&\left[ \tfrac{1}{\alpha_k}c^{(i+s-k)}_{\alpha-e_k} \,F_1^{\alpha_1}\cdots F_s^{\alpha_s}\, ,\,h\right]\\
&=&
\tfrac{\alpha_1}{\alpha_k}c^{(i+s-k)}_{\alpha-e_k} F_1^{\alpha_1-1}F_2^{\alpha_2}\cdots F_s^{\alpha_s}[F_1,h]
 + 
\tfrac{\alpha_2}{\alpha_k}c^{(i+s-k)}_{\alpha-e_k}  F_1^{\alpha_1}F_2^{\alpha_2-1}F_3^{\alpha_3}\cdots F_s^{\alpha_s}[F_2,h]\nonumber\\
&& +\cdots +
\tfrac{\alpha_s}{\alpha_k}c^{(i+s-k)}_{\alpha-e_k} F_1^{\alpha_1}\cdots F_{s-1}^{\alpha_{s-1}}F_s^{\alpha_s-1}[F_s,h] \nonumber\\
&=&
c^{(i+s-1)}_{\alpha-e_1} F_1^{\alpha_1-1}F_2^{\alpha_2}\cdots F_s^{\alpha_s}[F_1,h] 
+ c^{(i+s-2)}_{\alpha-e_2} F_1^{\alpha_1}F_2^{\alpha_2-1}F_3^{\alpha_3}\cdots F_s^{\alpha_s}[F_2,h]\nonumber\\
&& 
+\cdots+
c^{(i)}_{\alpha-e_s} F_1^{\alpha_1}\cdots F_{s-1}^{\alpha_{s-1}}F_s^{\alpha_s-1}[F_s,h] .\nonumber
\end{eqnarray}
Indeed, by assumption (C1), for $l\in\{ 1,\ldots,s \},$ we have
\begin{eqnarray}
[(\alpha_l-1)+1]c^{(i+s-k)}_{(\alpha-e_k-e_l)+e_l}=
[(\alpha_k-1)+1]c^{(i+s-l)}_{(\alpha-e_k-e_l)+e_k}.
\end{eqnarray}
If we have $\alpha_l=0$ in the above equality, then $(\alpha-e_k-e_l)+e_k\notin \mathbb{Z}\times\mathbb{N}^s$ and we can assume that $c^{(i+s-l)}_{(\alpha-e_k-e_l)+e_k}=0.$

Thus, the summand (\ref{row_pochIloczynu}) can be written as in (\ref{row_GiGi1Gi2}). 
Since the summands (\ref{row_summFkFl}) are equal to zero, the result follows.
$\rightBox$


\begin{lemma}\label{lem_P1}
Assume that $h$ is a homogeneous polynomial that is no power of any other polynomial of degree $<\deg h,$ polynomials $G_i,\ldots, G_{i+s-1},$ defined as in Lemma \ref{Lem_sumGiGi1Gi2}, satisfy the condition (C1). 
If polynomial $G_{i-1}$ is such that $[h^r,G_{i-1}]+[F_s,G_i]+\cdots+[F_1,G_{i+s-1}]=0,$ for some $r\in\mathbb{N}_+,$  then
\[
G_{i-1}=\sum_{\alpha=(\alpha_0,\ldots,\alpha_s)\in\mathbb{Z}\times\mathbb{N}^s}
c^{(i-1)}_{\alpha}F_1^{\alpha_1}\cdots F_s^{\alpha_s}h^{\alpha_0}
\]
with $c^{(i-1)}_{\alpha}\in \mathbb{C}$ such that for all $\alpha=(\alpha_0,\ldots,\alpha_s)\in\mathbb{Z}\times\mathbb{N}^s$ we have:\vspace{0.3cm}\newline
(P1) \hspace{1cm} if $\alpha_l>0$ with $l\in\{ 1,\ldots,s \},$ then 
$c^{(i-1)}_{\alpha-re_0}=\tfrac{\alpha_0}{r\alpha_l}c^{(i+s-l)}_{\alpha-e_l}.$
\end{lemma}
\dow 
Since $[h^r,G_{i-1}]=rh^{r-1}[h,G_{i-1}]=[h,rh^{r-1}G_{i-1}],$ it follows by Lemma \ref{Lem_sumGiGi1Gi2} and Corollary \ref{Corr_H_reduction} that 
\begin{eqnarray}\label{row_WzorNaGiMinus1}
G_{i-1}&=&\frac{1}{rh^{r-1}}\left[ \sum_{\alpha\in\mathbb{Z}\times\mathbb{N}_+\times\mathbb{N}^{s-1}}\tfrac{\alpha_0}{\alpha_1}
c^{(i+s-1)}_{\alpha-e_1}
F_1^{\alpha_1}\cdots F_s^{\alpha_s}h^{\alpha_0-1}
\right. \\
&&+
\sum_{\alpha\in\mathbb{Z}\times\{0\}\times\mathbb{N}_+\times\mathbb{N}^{s-2}}\tfrac{\alpha_0}{\alpha_2}
c^{(i+s-2)}_{\alpha-e_2}
F_2^{\alpha_2}\cdots F_s^{\alpha_s}h^{\alpha_0-1}\nonumber\\
&& \vdots \nonumber \\ 
&& \left. +
\sum_{\alpha\in\mathbb{Z}\times\{(0,\ldots,0)\}\times\mathbb{N}_+}\tfrac{\alpha_0}{\alpha_s}
c^{(i)}_{\alpha-e_s}
F_s^{\alpha_s}h^{\alpha_0-1} \right]+
P(h),\nonumber
\end{eqnarray}
for some polynomial $P$ in one variable, over $k.$
Thus, if $l=1$ (and $\alpha_1>0$), then the equation in (P1) is a direct consequence of (\ref{row_WzorNaGiMinus1}).

Assume that $l>1$ and $\alpha_l>0.$ If $\alpha_1=\ldots=\alpha_{l-1}=0,$ then the equality in (P1) is, as above, a~direct consequence of (\ref{row_WzorNaGiMinus1}). 
Thus, we can assume that there is $k<l$ such that $k\geq 1$ and $\alpha_k>0.$ 
We can also assume that we have already proved that 
\begin{equation}
c_{\alpha-re_0}^{(i-1)}=\tfrac{\alpha_0}{r\alpha_k}c_{\alpha-e_k}^{(i+s-k)}.
\end{equation} 
Now, using (C1) we obtain that
\begin{equation}
\left[ (\alpha_k-1)+1\right] c^{(i+s-l)}_{(\alpha-e_l-e_k)+e_k}
=
\left[ (\alpha_l-1)+1\right] c^{(i+s-k)}_{(\alpha-e_l-e_k)+e_l}
\end{equation}
and so
\begin{equation}
\tfrac{\alpha_0}{r\alpha_l} c^{(i+s-l)}_{\alpha-e_l}
=
\tfrac{\alpha_0}{r\alpha_k} c^{(i+s-k)}_{\alpha-e_k}.
\end{equation}
This completes the proof.
$\rightBox$

\begin{lemma}\label{lem_indukcjaC1C2} 
Assume that $r\in\mathbb{N}_+,$ $h$ is a homogeneous polynomial that is no power of any other polynomial of degree $<\deg h,$ polynomials $G_i,$ $\ldots,$ $G_{i+s-1},$ defined as in Lemma~\ref{Lem_sumGiGi1Gi2}, satisfy the condition (C1) and with polynomial $G_{i+s},$
defined similarly as in Lemma~\ref{Lem_sumGiGi1Gi2}, satisfy the following condition:\vspace{0.3cm}\newline
(C2) \hspace{1cm}\begin{minipage}{13cm} for all $1\leq l, k\leq s$ and $\alpha= (\alpha_0,\ldots,\alpha_s)\in\mathbb{Z}\times\mathbb{N}^s$ \newline if $\alpha_k>0,$ then $c^{(i+l-1)}_{\alpha-re_0}= \tfrac{\alpha_0}{r\alpha_k}c^{(i+l+s-k)}_{\alpha-e_k}$ 
\end{minipage}
 \vspace{0.3cm}\newline
 If polynomial $G_{i-1}$ is such that $[h^r,G_{i-1}]+[F_s,G_i]+\ldots+[F_1,G_{i+s-1}]=0,$ then\newline
(a) $G_{i-1},G_i,\ldots,G_{i+s-2}$ satisfy the condition (C1),\newline
(b) $G_{i-1},G_i,\ldots,G_{i+s-1}$  satisfy the condition (C2),\newline
all with the same $r$ and $i$ replaced by $i-1.$
\end{lemma}
\dow We start with the proof of (b). Take any $1\leq l, k\leq s$ and $\alpha=(\alpha_0,\ldots,\alpha_s)\in\mathbb{Z}\times\mathbb{N}^s$ with $\alpha_k>0.$ 

Assume that $l\geq 2.$ Then, by (C2) we have
\begin{eqnarray}
c^{((i-1)+l-1)}_{\alpha-re_0}
&=&c^{(i+(l-1)-1)}_{\alpha-re_0}
=\tfrac{\alpha_0}{r\alpha_k}c^{(i+(l-1)+s-k)}_{\alpha -e_k}
=\tfrac{\alpha_0}{r\alpha_k}c^{((i-1)+l+s-k)}_{\alpha -e_k}.
\end{eqnarray}

For $l=1$ we can use property (P1) with $l$ replaced by $k.$

Now, we prove (a). Take any $1\leq k,l\leq s$ and $\alpha=(\alpha_0,\ldots,\alpha_s)\in\mathbb{Z}\times\mathbb{N}^s.$
Let us notice that
\begin{eqnarray}\label{row_ProofC1_1}
(\alpha_l+1)c^{((i-1)+s-k)}_{\alpha+e_l}
=(\alpha_l+1)c^{(i+(s-k)-1)}_{\alpha+e_l}
=(\alpha_l+1)\tfrac{\alpha_0}{r(\alpha_l+1)}c^{(i+(s-k)+s-l)}_{\alpha+re_0}.
\end{eqnarray}
Indeed, if $k=s,$ then it is a consequence of (P1):
\begin{eqnarray}
(\alpha_l+1)c^{(i+(s-k)-1)}_{\alpha+e_l}
&=&
(\alpha_l+1)c^{(i-1)}_{\alpha+e_l}
= (\alpha_l+1)c^{(i-1)}_{(\alpha+e_l+re_0)-re_0}
\\
\nonumber
&=&
 (\alpha_l+1)\tfrac{\alpha_0}{r(\alpha_l+1)}c^{(i+s-l)}_{\alpha+re_0}
=
 (\alpha_l+1)\tfrac{\alpha_0}{r(\alpha_l+1)}c^{(i+(s-k)+s-l)}_{\alpha+re_0}
\end{eqnarray} 
and if $k<s$ it follows from (C2):

\begin{eqnarray}
(\alpha_l+1)c^{(i+(s-k)-1)}_{\alpha+e_l}
&=&
(\alpha_l+1)c^{(i+(s-k)-1)}_{(\alpha+e_l+re_0)-re_0}
=
(\alpha_l+1)\tfrac{\alpha_0}{r(\alpha_l+1)}c^{(i+(s-k)+s-l)}_{\alpha+re_0}.
\end{eqnarray}
Similarly, we have
\begin{eqnarray}\label{row_ProofC1_2}
(\alpha_k+1)c^{((i-1)+s-l)}_{\alpha+e_k}
=(\alpha_k+1)c^{(i+(s-l)-1)}_{\alpha+e_k}
=(\alpha_k+1)\tfrac{\alpha_0}{r(\alpha_k+1)}c^{(i+(s-l)+s-k)}_{\alpha+re_0}.
\end{eqnarray}
By (\ref{row_ProofC1_1}) and (\ref{row_ProofC1_2}), we have
\begin{eqnarray}
(\alpha_l+1)c^{((i-1)+s-k)}_{\alpha+e_l}=(\alpha_k+1)c^{((i-1)+s-l)}_{\alpha+e_k}.
\end{eqnarray} 
$\rightBox$

The following lemma is obvious.

\begin{lemma}\label{lem_StartIndukcji}
Assume that $h$ is a homogeneous polynomial such that $\deg h$ divides $N.$
Let $G_{N}=\alpha_{N}h^{\frac{N}{\deg h}},$ $G_{N+1}=\ldots=G_{N+s}=0.$  Then\newline
(a) polynomials $G_{N},G_{N+1},\ldots,G_{N+s-1}$  satisfy the condition (C1),\newline
(b) for any $r\in\mathbb{N}_+,$ polynomials  $G_{N},G_{N+1},\ldots,G_{N+s}$  satisfy the condition (C2),\newline
all with $i$ replaced by $N.$
\end{lemma}



\section{Formulas for $G_j$}\label{Sect_formula_for_Gj}
The main result of this section is the following formula for polynomials $G_j.$

\begin{theorem}\label{prop_formula}
Let $F=F_1+\cdots+F_{s}+F_d$ and $G=G_1+ \cdots +G_N$ be
elements of $\mathbb{C}[x_1,\ldots,x_n],$ where $F_i$ and $G_j$ are homogeneous polynomials of degree $i$ and $j,$ respectively.
Assume that 
\begin{equation}
F_d=h^d \qquad\text{and}\qquad G_N=a_Nh^N
\end{equation}
for some homogeneous polynomial $h$ of degree 1 and $a_N\in \mathbb{C}\setminus\{ 0\} .$

If $\deg [F,G]<d+i$ for $i\in\{1,\ldots,N \} ,$ then there are $a_i,\ldots,a_{N-1}\in \mathbb{C}$ such that for $j=i,i+1,\ldots,N,$ 
we have 
\begin{equation}\label{row_formula1}
G_j=\sum_{\substack{\alpha=(\alpha_0,\ldots,\alpha_s)\in\mathbb{Z}\times\mathbb{N}^s\\
\alpha_0+\alpha_1+2\alpha_2+\cdots+s\alpha_s=j\\
\alpha_0+d(\alpha_1+\cdots+\alpha_s)\leq N}}
c^{(j)}_{\alpha}F_1^{\alpha_1}\cdots F_s^{\alpha_s}h^{\alpha_0},
\end{equation}
where for $\alpha=(\alpha_0,\ldots,\alpha_s)\in\mathbb{Z}\times\mathbb{N}^s$ with $\alpha_0+\alpha_1+2\alpha_2+\cdots+s\alpha_s=j$ and $\alpha_0+d(\alpha_1+\cdots+\alpha_s)\leq N,$ we have
\begin{equation}\label{row_formula2}
c^{(j)}_{\alpha}=
\frac{\prod_{k=1}^{\alpha_1+\cdots+\alpha_s}(\alpha_0+dk)}
{\left( \prod_{i_1=1}^{\alpha_1} di_1\right) \cdots 
\left( \prod_{i_s=1}^{\alpha_s} di_s\right)} a_{\alpha_0+d(\alpha_1+\cdots+\alpha_s)}.
\end{equation}
\end{theorem}
\dow
In what follows, we assume that $c^{(j)}_{(\alpha_0,\ldots,\alpha_s)}=0$ if $\alpha_0+\alpha_1+2\alpha_2+\cdots+s\alpha_s\neq j$ or $\alpha_0+d(\alpha_1+\cdots+\alpha_s)> N.$

Let us notice that the above formula holds for $G_N.$ 
Indeed, if $(\alpha_1,\ldots,\alpha_s)\in\mathbb{N}^s\setminus\{ (0,\ldots,0) \}$ and $\alpha_0+\alpha_1+2\alpha_2+\cdots+s\alpha_s=N,$  then $\alpha_0+d(\alpha_1+\cdots+\alpha_s)>N.$ 
Hence, in the sum (\ref{row_formula1}), for $G_N,$ there is only one summand $a_Nh^N.$
Let us also notice that the above formula holds for the polynomials $G_{N+1}=\ldots=G_{N+s}=0.$

By Lemma~\ref{lem_StartIndukcji}, we know that the polynomials $G_N,\ldots,G_{N+s-1}$ satisfy the condition (C1), and that the polynomials $G_N,\ldots,G_{N+s}$ satisfy the condition (C2).
Assume that $i\leq N-1.$ Then, $\deg [F,G]<d+N-1$ and so $[h^d,G_{N-1}]+[F_s,G_N]+\cdots+[F_1,G_{N+s-1}]=0.$
Since $G_N,\ldots,G_{N+s-1}$ satisfy the condition (C1), we see from Lemma~\ref{lem_P1} that the polynomial $G_{N-1}$ is of the form
\[
G_{N-1}=\sum_{\substack{\alpha=(\alpha_0,\ldots,\alpha_s)\in\mathbb{Z}\times\mathbb{N}^s\\
\alpha_0+\alpha_1+2\alpha_2+\cdots+s\alpha_s=N-1}}
c^{(N-1)}_{\alpha}F_1^{\alpha_1}\cdots F_s^{\alpha_s}h^{\alpha_0}
\]
and satisfy (P1).
But, since $G_N=a_Nh^N$ and $G_{N+1}=\ldots=G_{N+s-1}=0,$ we see, from (P1), that if $(\alpha_1,\ldots,\alpha_s)\in\mathbb{N}^s\setminus\{ (0,\ldots,0),(0,\ldots,0,1) \},$ then $c^{(N-1)}_{(\alpha_0,\ldots,\alpha_s)}=0.$
For $\alpha_s=1,$ by (P1), we have
\[
c^{(N-1)}_{(N-d,0,\ldots,0,1)}=\tfrac{N}{d\cdot 1}c^{(N)}_{(N,0,\ldots,0)}=\tfrac{N}{d\cdot 1}a_N=\tfrac{(N-d)+d\cdot 1}{d\cdot 1}a_{(N-d)+d(0+\cdots+0+1)}.
\]
Setting  $a_{N-1}=c^{(N-1)}_{(N-1,0,\ldots,0)},$ we obtain that $G_{N-1}$ satisfy formulas (\ref{row_formula1}) and (\ref{row_formula2}).
By Lemma~\ref{lem_indukcjaC1C2}, we also know that the polynomials $G_{N-1},\ldots,G_{N+s-2}$ satify the condition (C1), and that the polynomials $G_{N-1},\ldots,G_{N+s-1}$ satify the condition (C2). Thus, we can use the induction.

Now, assume that for some $l\in\{ i+1,i+2,\ldots,N \} ,$ polynomials $G_l,\ldots,G_{l+s}$ are given by the formulas (\ref{row_formula1}) and (\ref{row_formula2}). 
Assume  also that $G_l,\ldots,G_{l+s-1}$ satisfy the condition (C1) and $G_l,\ldots,G_{l+s}$ satisfy the condition (C2). 
Since $\deg [F,G]<d+i\leq d+(l-1),$ it follows that $[h^d,G_{l-1}]+[F_s,G_l]+\cdots+[F_1,G_{l+s-1}]=0.$
Hence, by Lemma~\ref{lem_P1}, the polynomial $G_{l-1}$ is of the form
\begin{equation}\label{row_PolyLMinus1}
G_{l-1}=\sum_{\substack{\alpha=(\alpha_0,\ldots,\alpha_s)\in\mathbb{Z}\times\mathbb{N}^s\\
\alpha_0+\alpha_1+2\alpha_2+\cdots+s\alpha_s=l-1}}
c^{(l-1)}_{\alpha}F_1^{\alpha_1}\cdots F_s^{\alpha_s}h^{\alpha_0}
\end{equation}
and satisfy (P1). 
Moreover, by Lemma~\ref{lem_indukcjaC1C2}, the polynomials  $G_{l-1},\ldots,G_{l+s-2}$ satisfy the condition (C1) and $G_{l-1},\ldots,G_{l+s-1}$ satisfy the condition (C2).
Since $G_{l-1}$ satisfy (P1), we have that  
\begin{eqnarray}\label{row_GLminus1Formula}
G_{l-1}&=&
a_{l-1}h^{l-1}+
\sum_{\substack{\alpha=(\alpha_0,\ldots,\alpha_s)\in\mathbb{Z}\times\mathbb{N}_+\times\mathbb{N}^{s-1}\\
\alpha_0+\alpha_1+2\alpha_2+\cdots+s\alpha_s=l-1}}
\tfrac{\alpha_0+d}{d\alpha_1} c^{(l+s-1)}_{(\alpha_0+d,\alpha_1-1,\alpha_2,\ldots,\alpha_s)}F_1^{\alpha_1}\cdots F_s^{\alpha_s}h^{\alpha_0}\\
&&
+
\sum_{\substack{\alpha=(\alpha_0,\ldots,\alpha_s)\in\mathbb{Z}\times\{ 0 \}\times\mathbb{N}_+\times\mathbb{N}^{s-2}\\
\alpha_0+\alpha_1+2\alpha_2+\cdots+s\alpha_s=l-1}}
\tfrac{\alpha_0+d}{d\alpha_2} 
c^{(l+s-2)}_{(\alpha_0+d,0,\alpha_2-1,\alpha_3,\ldots,\alpha_s)}F_2^{\alpha_2}\cdots F_s^{\alpha_s}h^{\alpha_0}
\nonumber\\
&&
\vdots
\nonumber\\
&&
+
\sum_{\substack{\alpha=(\alpha_0,\ldots,\alpha_s)\in\mathbb{Z}\times\{ (0,\ldots,0) \}\times\mathbb{N}_+\\
\alpha_0+\alpha_1+2\alpha_2+\cdots+s\alpha_s=l-1}}
\tfrac{\alpha_0+d}{d\alpha_s} 
c^{(l+s-s)}_{(\alpha_0+d,0,\ldots,0,\alpha_s-1)}F_s^{\alpha_s}h^{\alpha_0}.
\nonumber
\end{eqnarray}
Notice that $c^{(l+s-1)}_{(\alpha_0+d,\alpha_1-1,\alpha_2,\ldots,\alpha_s)}$ can be different from zero only if $(\alpha_0+d)+d(\alpha_1-1+\alpha_2+\cdots+\alpha_s)\leq N$ 
(because there is no components with degree $>N$ in the polynomial $G$),
that is only if $\alpha_0+d(\alpha_1+\alpha_2+\cdots+\alpha_s)\leq N.$ 
The same holds for $c^{(l+s-2)}_{(\alpha_0+d,0,\alpha_2-1,\alpha_3\,\ldots,\alpha_s)},$ $\ldots,$  $c^{(l)}_{(\alpha_0+d,0,\ldots,0,\alpha_s-1)}.$ 
Hence, in (\ref{row_PolyLMinus1}), we can add the restriction $\alpha_0+d(\alpha_1+\cdots+\alpha_s)\leq N.$

Assume that $\alpha_1>0,$ $ \alpha_0+\alpha_1+2\alpha_2+\cdots+s\alpha_s=l-1$ and $\alpha_0+d(\alpha_1+\cdots+\alpha_s)\leq N.$ 
Then, by (\ref{row_GLminus1Formula}) and the above assumptions, we have
\begin{eqnarray}
c^{(l-1)}_{(\alpha_0,\ldots,\alpha_s)}
&=&
\tfrac{\alpha_0+d}{d\alpha_1}c^{(l+s-1)}_{(\alpha_0+d,\alpha_1-1,\alpha_2,\ldots,\alpha_s)}
\\
&=&
\tfrac{\alpha_0+d}{d\alpha_1}
\frac{\prod_{k=1}^{\alpha_1-1+\cdots+\alpha_s}((\alpha_0+d)+dk)}
{\left( \prod_{i_1=1}^{\alpha_1-1}di_1\right)\cdots
\left( \prod_{i_s=1}^{\alpha_s}di_s\right)} 
a_{(\alpha_0+d)+d((\alpha_1-1)+\alpha_2+\cdots+\alpha_s)}
\nonumber
\\
&=&
\frac{\prod_{k=1}^{\alpha_1+\cdots+\alpha_s}(\alpha_0+dk)}
{\left( \prod_{i_1=1}^{\alpha_1}di_1\right)\cdots
\left( \prod_{i_s=1}^{\alpha_s}di_s\right)} 
a_{\alpha_0+d(\alpha_1+\alpha_2+\cdots+\alpha_s).}
\nonumber
\end{eqnarray}

Now, assume that $\alpha_1=0,$ $\alpha_2>0,$  $ \alpha_0+2\alpha_2+\cdots+s\alpha_s=l-1$ and $\alpha_0+d(\alpha_2+\cdots+\alpha_s)\leq N.$
Similarly, as before, we have
\begin{eqnarray}
c^{(l+s-2)}_{(\alpha_0,\ldots,\alpha_s)}
&=&
\tfrac{\alpha_0+d}{d\alpha_2}c^{(l+1)}_{(\alpha_0+d,0,\alpha_2-1,\alpha_3,\ldots,\alpha_s)}
\\
&=&
\tfrac{\alpha_0+d}{d\alpha_2}
\frac{\prod_{k=1}^{\alpha_2-1+\alpha_3+\cdots+\alpha_s}((\alpha_0+d)+dk)}
{\left( \prod_{i_2=1}^{\alpha_2-1}di_2\right) \cdots 
\left( \prod_{i_s=1}^{\alpha_s}di_s \right)} a_{(\alpha_0+d)+d((\alpha_2-1)+\alpha_3+\cdots+\alpha_s)}
\nonumber
\\
&=&
\frac{\prod_{k=1}^{\alpha_2+\alpha_3+\cdots+\alpha_s}(\alpha_0+dk)}
{\left( \prod_{i_2=1}^{\alpha_2}di_2\right) \cdots 
\left( \prod_{i_s=1}^{\alpha_s}di_s \right)} a_{\alpha_0+d(\alpha_2+\alpha_3+\cdots+\alpha_s)}
\nonumber
\end{eqnarray}

The similar calculation shows that formula (\ref{row_formula2}) holds also in the cases $\alpha_1=\ldots=\alpha_q=0,$ and $\alpha_{q+1}>0,$  with $\alpha_0+\alpha_1+2\alpha_2+\cdots+s\alpha_s=l-1$ and $\alpha_0+d(\alpha_1+\cdots+\alpha_s)\leq N,$ for all $q=3,\ldots,s-1.$

Since, by Lemma~\ref{lem_indukcjaC1C2} we can use induction, the result follows. $\rightBox$

Using similar arguments as above one can prove the following result in which we consider the case $\deg h\geq 1;$ interested one is $\deg h>1.$
Of course, the following theorem is a generalization of the above one.

\begin{theorem}\label{prop_formulaB}
Let $F=F_1+\cdots+F_{s}+F_d$ and $G=G_1+\cdots +G_N$ be
elements of $\mathbb{C}[x_1,\ldots,x_n],$ where $F_i$ and $G_j$ are homogeneous polynomials of degree $i$ and $j,$ respectively.
Assume that $h\in \mathbb{C}[x_1,\ldots,x_n]$ is a homogeneous polynomial that is no power of any other polynomial of degree $<\deg h,$  such that $\deg h| \gcd(d,N)$ and 
\begin{equation}
F_d=h^{d/\deg h} \qquad\text{and}\qquad G_N=a_Nh^{N/\deg h}
\end{equation}
for some $a_N\in \mathbb{C}\setminus\{ 0\} .$
In what follows, we put $r=\tfrac{d}{\deg h}$ and $t=\deg h.$

If $\deg [F,G]<d+i$ for $i\in\{1,\ldots,N \} ,$ then there are $a_i,\ldots,a_{N-1}\in \mathbb{C}$ such that for $j=i,i+1,\ldots,N,$ 
we have 
\begin{equation}\label{row_formula1B}
G_j=\sum_{\substack{\alpha=(\alpha_0,\ldots,\alpha_s)\in\mathbb{Z}\times\mathbb{N}^s\\
t\alpha_0+\alpha_1+2\alpha_2+\cdots+s\alpha_s=j\\
t\alpha_0+d(\alpha_1+\cdots+\alpha_s)\leq N}}
c^{(j)}_{\alpha}F_1^{\alpha_1}\cdots F_s^{\alpha_s}h^{\alpha_0},
\end{equation}
where for $\alpha=(\alpha_0,\ldots,\alpha_s)\in\mathbb{Z}\times\mathbb{N}^s$ with $t\alpha_0+\alpha_1+2\alpha_2+\cdots+s\alpha_s=j$ and $t\alpha_0+d(\alpha_1+\cdots+\alpha_s)\leq N,$ 
we have
\begin{equation}\label{row_formula2B}
c^{(j)}_{\alpha}=
\frac{\prod_{k=1}^{\alpha_1+\cdots+\alpha_s}(\alpha_0+rk)}
{\left( \prod_{i_1=1}^{\alpha_1} ri_1\right) \cdots 
\left( \prod_{i_s=1}^{\alpha_s} ri_s\right)} a_{t\alpha_0+d(\alpha_1+\cdots+\alpha_s)}.
\end{equation}

Moreover, we have $a_j=0$ when $\deg h\nmid j.$ 
\end{theorem}

\section{Divisibility of $F_s$}\label{sect_divisibility_Fs}

In this section, we show that, if $\deg [F,G]$ is small enough then $F_{s}$ is divisible by $h$ (see Theorem~\ref{prop_F3equal_f2h}) or even by suitable power of $h$ (see Theorem~\ref{prop_F3equal_f2h2}).
This will be established in the case $\deg h=1$ and in the case $\deg h>1$ with the additional assumption that $h$ is a square-free polynomial (see Theorem~\ref{prop_divisibility_Fs_h2_t}).
In the case $\deg h>1$ but without assumption that $h$ is square-free, we obtain that $F_s$ is divisible by $\sqrf h,$ 
where $\sqrf h$ is defined as the product of all irreducible factors of the polynomial $h$ but taken without multiplicity, for example $\sqrf\left( x_1^3(x_2-1)^2\right)=x_1(x_2-1).$
A polynomial $\sqrf h$ for a given $h\in \mathbb{C}[x_1,\ldots,x_n]$ is defined only up to a constant factor $u\in \mathbb{C}\setminus\{0\},$ but for our purpose it is not a disadvantage.
If the polynomial $h$ is itself a \textit{square-free} polynomial, that is with irreducible factors with multiplicity one, then we put $\sqrf h=h.$
For a general definition of square-free elements of any commutative rings we refer the reader to \cite{Jedrzejewicz_Zielinski}.

\subsection{The case $\deg h=1$}
In this subsection we consider the case $\deg h=1.$
First of all, let us notice that if $\deg h\leq 2$ and $h$ is no power of any other polynomial with degree $<\deg h,$ then the polynomial $h$ is square-free.

\begin{theorem}\label{prop_F3equal_f2h}
Let $F$ and $G$ be
as in Theorem \ref{prop_formula}.
If $N\not\equiv 0\, (\mod d)$ and $\deg[F,G]< d+i$ with $i<\tfrac{s}{d}N=\tfrac{d-1}{d}N,$ then $F_s$ is divisible by $h.$ In other words, there is a homogeneous polynomial $\tilde{F}_{s-1}$ of degree $s-1$ such that $F_s=\tilde{F}_{s-1}h.$
\end{theorem}

In the proof of the above theorem we use the following

\begin{lemma}\label{lem_MinimalValueAlpha0}
For all integers $1\leq j< N,$ the minimal value
\begin{equation}
\min \{\, \alpha_0\,|\, (\alpha_0,\ldots,\alpha_s)\in
\mathbb{Z}\times\mathbb{N}^s,\, \alpha_0+\alpha_1+2\alpha_2+\cdots+s\alpha_s=j,\,
\alpha_0+d(\alpha_1+\cdots+\alpha_s)\leq N\,\}
\end{equation}
is equal to $dj-sN$ and is attained at only one point $(dj-sN,0,\ldots,0,N-j).$
\end{lemma}
\dow
Let $Z_{j,N,d}:=\{\,  (\alpha_0,\ldots,\alpha_s)\in
\mathbb{Z}\times\mathbb{N}^s\,|\, \alpha_0+\alpha_1+2\alpha_2+\cdots+s\alpha_s=j,\,
\alpha_0+d(\alpha_1+\cdots+\alpha_s)\leq N\,\} .$
Notice that for all $(\alpha_0,\ldots,\alpha_s)\in Z_{j,N,d},$ we have $s\alpha_1+(s-1)\alpha_2+\cdots+2\alpha_{s-1}+\alpha_s\leq N-j$ 
(in particular, $\alpha_s\leq N-j$).
Thus, $s\alpha_s\leq s(N-j) -s\left(s\alpha_1+(s-1)\alpha_2+\cdots+2\alpha_{s-1}\right)$ 
and $\alpha_1+2\alpha_2+\cdots+s\alpha_s\leq s(N-j) -[s^2-1]\alpha_1 -[s(s-1)-2]\alpha_2-\cdots-[s\cdot 2-(s-1)]\alpha_{s-1}.$
Since $s^2-1,s(s-1)-2,\ldots, s\cdot 2-(s-1)$ are strictly positive numbers, it follows that for $(\alpha_1,\ldots,\alpha_s)\in Z_{j,N,d}$ if $(\alpha_1,\ldots,\alpha_{s-1})\neq (0,\ldots,0)\in\mathbb{N}^{s-1},$ then $\alpha_0>dj-sN.$
Of course, if $(\alpha_0,0,\ldots,0,\alpha_s)\in Z_{j,N,d}$ with $\alpha_s<N-j,$ then also $\alpha_0=j-s\alpha_s>j-s(N-j)=dj-sN.$
$\rightBox$

Now we can prove Theorem~\ref{prop_F3equal_f2h}.\vspace{0.2cm}

\dow[ of Theorem \ref{prop_F3equal_f2h}]
Since $\deg [F,G]<d+i,$ we see by Theorem~\ref{prop_formula} that
\begin{equation}\label{row_propHDzieliF3}
G_i=\sum_{\substack{\alpha=(\alpha_0,\ldots,\alpha_s)\in\mathbb{Z}\times\mathbb{N}^s\\
\alpha_0+\alpha_1+2\alpha_2+\cdots+s\alpha_s=i\\
\alpha_0+d(\alpha_1+\cdots+\alpha_s)\leq N}}
c^{(i)}_{\alpha}F_1^{\alpha_1}\cdots F_s^{\alpha_s}h^{\alpha_0}.
\end{equation}
By Lemma~\ref{lem_MinimalValueAlpha0}, all summands of (\ref{row_propHDzieliF3}), except $c^{(i)}_{(di-sN,0,\ldots,0,N-i)}F_s^{N-i}h^{di-sN},$ involves $h^u$ with $u>di-sN.$
Let us also notice that $di-sN<0.$
Thus, multiplying both sides of (\ref{row_propHDzieliF3}) by $h^{sN-di},$ we obtain that $h$ divides $c^{(i)}_{(di-sN,0,\ldots,0,N-i)}F_s^{N-i}.$ 
Since $N\not\equiv  0 (\mod d),$ we see that $-sN\not\equiv 0 (\mod d)$ and so $c^{(i)}_{(di-sN,0,\ldots,0,N-i)}=\frac{\prod_{k=1}^{N-i}(di-sN+dk)}{\prod_{i_s=1}^{N-i}di_s}a_N\neq 0.$
Thus, we see that $h|F_s^{N-i}$ and so $h|F_s$ because $h$ is a square-free polynomial.
$\rightBox$ 

By Theorems~\ref{prop_formula} and \ref{prop_F3equal_f2h}, we obtain the following 

\begin{corollary}
Let $F$ and $G$ be
as in Theorem \ref{prop_formula}.
If $N\not\equiv  0\, (\mod d)$ and $\deg[F,G]< d+i$ with $i<\tfrac{s}{d}N,$
then for $j=i,i+1,\ldots,N,$ we have 
\begin{equation}
G_j=\sum_{\substack{\alpha=(\alpha_0,\ldots,\alpha_s)\in\mathbb{Z}\times\mathbb{N}^s\\
\alpha_0+\alpha_1+2\alpha_2+\cdots+s\alpha_s=j\\
\alpha_0+d(\alpha_1+\cdots+\alpha_s)\leq N}}
c^{(j)}_{\alpha}F_1^{\alpha_1}\cdots{F}_{s-1}^{\alpha_{s-1}}\tilde{F}_{s-1}^{\alpha_s}h^{\alpha_0+\alpha_s},
\end{equation}
where $a_i,\ldots,a_{N-1}\in \mathbb{C}$ and $c^{(j)}_{\alpha}$ are as in Theorem \ref{prop_formula}, and $\tilde{F}_{s-1}$ is as in Theorem~\ref{prop_F3equal_f2h}.
\end{corollary}

The following lemma we use in the proof of Theorem~\ref{prop_F3equal_f2h2} below.

\begin{lemma}
For all $k\in\mathbb{N}$ such that $d>2k$ and for all integers $1\leq j< N,$ the minimal value
\begin{equation}
\min \{\, \alpha_0+k\alpha_s\,|\, (\alpha_0,\ldots,\alpha_s)\in Z_{j,N,d}\,\}
\end{equation}
is equal to $(d-k)j-(s-k)N$ and is attained at only one point $(dj-sN,0,\ldots,0,N-j).$
\end{lemma}
\dow
As in the proof of Lemma~\ref{lem_MinimalValueAlpha0}, we have 
$$
(s-k)\alpha_s\leq (s-k)(N-j)-(s-k)\left( s\alpha_1+(s-1)\alpha_2+\cdots+2\alpha_{s-1}\right).
$$
From this we see that
$
\alpha_1+2\alpha_2+\cdots+(s-1)\alpha_{s-1}+(s-k)\alpha_s \leq 
(s-k)(N-j)- [(s-k)s-1]\alpha_1 -[(s-k)(s-1)-2]\alpha_2-\cdots -[(s-k)\cdot 2-(s-1)]\alpha_{s-1}
$
and so
\begin{eqnarray}
\alpha_0+k\alpha_s
&=&
j-[\alpha_1+2\alpha_2+\cdots+(s-1)\alpha_{s-1}+(s-k)\alpha_s]
\\
\nonumber
&\geq&
j-(s-k)(N-j)+[(s-1)s-1]\alpha_1 
\\
\nonumber
&&
+[(s-k)(s-1)-2]\alpha_2+\cdots +[(s-k)\cdot 2-(s-1)]\alpha_{s-1}
\end{eqnarray}
Since $(s-k)s-1>(s-k)(s-1)-2>\ldots>(s-k)\cdot 2-(s-1)=s-(2k-1)>0,$ it follows that 
 for $(\alpha_0,\ldots,\alpha_s)\in Z_{j,N,d}$ if $(\alpha_1,\ldots,\alpha_{s-1})\neq (0,\ldots,0)\in\mathbb{N}^{s-1},$ then $\alpha_0+k\alpha_s>(d-1)j-(s-1)N.$
Of course, if $(\alpha_0,0,\ldots,0,\alpha_s)\in Z_{j,N,d}$ with $\alpha_s<N-j,$ then also $\alpha_0+k\alpha_s=j-(s-k)\alpha_s>j-(s-k)(N-j)=(d-k)j-(s-k)N.$
$\rightBox$

By the above lemma, one can use, several times, the similar arguments as in the proof of Theorem~\ref{prop_F3equal_f2h} to show the following result (see also proof of Theorem~\ref{prop_divisibility_Fs_h2_t} in the next subsection).

\begin{theorem}\label{prop_F3equal_f2h2}
Let $F,G$ and $h$ be
as in Theorem \ref{prop_formula}.
If $N\not\equiv 0\, (\mod d)$ and $\deg[F,G]< d+i$ with $i<\tfrac{s-k}{d-k}N$ with $2k<d,$ then $F_s$ is divisible by $h^{k+1}.$ In other words, there is a homogeneous polynomial $\tilde{F}_{s-k-1}$ of degree $s-k-1$ such that $F_s=\tilde{F}_{s-k-1}h^{k+1}.$
\end{theorem}

By Theorems~\ref{prop_formula} and \ref{prop_F3equal_f2h2}, we obtain the following 

\begin{corollary}
Let $F$ and $G$ be
as in Theorem \ref{prop_formula}.
If $N\not\equiv 0\, (\mod d)$ and $\deg[F,G]< d+i$ with $i<\tfrac{s-k}{d-k}N$ with $2k<d,$
then we have 
\begin{equation}
G_j=\sum_{\substack{\alpha=(\alpha_0,\ldots,\alpha_s)\in\mathbb{Z}\times\mathbb{N}^s\\
\alpha_0+\alpha_1+2\alpha_2+\cdots+s\alpha_s=j\\
\alpha_0+d(\alpha_1+\cdots+\alpha_s)\leq N}}
c^{(j)}_{\alpha}F_1^{\alpha_1}\cdots F_{s-1}^{\alpha_{s-1}}\tilde{F}_{s-k-1}^{\alpha_s}h^{\alpha_0+(k+1)\alpha_s},
\end{equation}
where $a_i,\ldots,a_{N-1}\in \mathbb{C}$ and $c^{(j)}_{\alpha}$ are as in Theorem \ref{prop_formula} and $\tilde{F}_{s-k-1}$ is as in Theorem~\ref{prop_F3equal_f2h2}.
\end{corollary}

\subsection{The case $\deg h>1$}
The results given in this subsection are also true for $\deg h=1,$ but of course are interested in the case $\deg h>1.$
We start with the following lemma.

\begin{lemma}\label{lem_MinimalValueAlpha0Plus_kAlpha_s_t}
Let $d,t$ and $N$ be as in Theorem~\ref{prop_formulaB}.
For all $k\in\mathbb{N}$ such that $d>2kt$ and for all integers $1\leq j< N,$ the minimal value
\begin{equation}
\min \{\, \alpha_0+k\alpha_s\,|\, (\alpha_0,\ldots,\alpha_s)\in Z_{j,N,d}^{(t)}\,\}
\end{equation}
is equal to $\tfrac{1}{t}\left[ (d-kt)j-(s-kt)N\right]$ and is attained at only one point $(\tfrac{1}{t}(dj-sN),0,\ldots,0,N-j),$ where
$
Z_{j,N,d}^{(t)}=\left\{ (\alpha_0,\ldots,\alpha_s)\in\mathbb{Z}\times\mathbb{N}^s \,:\, t\alpha_0+\alpha_1+2\alpha_2+\cdots+s\alpha_s=j, t\alpha_0+d(\alpha_1+\cdots+\alpha_s)\leq N \right\}
$
\end{lemma}
\dow
As in the proof of Lemma~\ref{lem_MinimalValueAlpha0}, we have 
$$
(s-kt)\alpha_s\leq (s-kt)(N-j)-(s-kt)\left( s\alpha_1+(s-1)\alpha_2+\cdots+2\alpha_{s-1}\right).
$$
From this we see that
$
\alpha_1+2\alpha_2+\cdots+(s-1)\alpha_{s-1}+(s-kt)\alpha_s \leq 
(s-kt)(N-j)- [(s-kt)s-1]\alpha_1 -[(s-kt)(s-1)-2]\alpha_2-\cdots -[(s-kt)\cdot 2-(s-1)]\alpha_{s-1}.
$

Since $\alpha_0=\tfrac{1}{t}\left[j-(\alpha_1+2\alpha_2+\cdots+s\alpha_s) \right],$ it follows that
\begin{eqnarray}
\alpha_0+k\alpha_s
&=&
\tfrac{1}{t}\left[j-(\alpha_1+2\alpha_2+\cdots+s\alpha_s) \right]+\frac{kt\alpha_s}{t}
\\
\nonumber
&=&
\tfrac{1}{t}\left[j-(\alpha_1+2\alpha_2+\cdots+(s-1)\alpha_{s-1}+(s-kt)\alpha_s) \right]
\\
\nonumber
&\geq&
\tfrac{1}{t}\left\{
j+(s-kt)(j-N)+[(s-kt)s-1]\alpha_1+[(s-kt)(s-1)-2]\alpha_2+\right.
\\
\nonumber
&&
\left. \cdots+[(s-kt)2-(s-1)]\alpha_{s-1}
\right\}.
\end{eqnarray}
Since $(s-kt)s-1>(s-kt)(s-1)-2>\ldots>(s-kt)\cdot 2-(s-1)=s-(2kt-1)>0,$ it follows that 
 for $(\alpha_0,\ldots,\alpha_s)\in Z_{j,N,d}^{(t)}$ if $(\alpha_1,\ldots,\alpha_{s-1})\neq (0,\ldots,0)\in\mathbb{N}^{s-1},$ then $\alpha_0+k\alpha_s>\tfrac{1}{t}\left\{(d-kt)j-(s-kt)N\right\}.$
Of course, if $(\alpha_0,0,\ldots,0,\alpha_s)\in Z_{j,N,d}^{(t)}$ with $\alpha_s<N-j,$ then also $\alpha_0+k\alpha_s=\tfrac{1}{t}\left\{ j-(s-kt)\alpha_s\right\}>\tfrac{1}{t}\left\{j-(s-k)(N-j\right\})=\tfrac{1}{t}\left\{(d-kt)j-(s-kt)N\right\}.$
$\rightBox$
\newline

Now, we can use the above lemma and Theorem~\ref{prop_formulaB} to prove the following result.

\begin{theorem}\label{prop_divisibility_Fs_h2_t}
Let $F,G$ and $h$ be as in Theorem \ref{prop_formulaB}.
Assume also that the polynomial $h$ is square-free.
If $N\not\equiv  0\, (\mod d)$ and $\deg[F,G]< d+i$ with $i<\tfrac{s-k}{d-k}N$ and $2kt<d,$ then $F_s$ is divisible by $h^{k+1}.$ In other words, there is a homogeneous polynomial $\tilde{F}_{s-t(k+1)}$ of degree $s-t(k+1)$ such that $F_s=\tilde{F}_{s-t(k+1)}h^{k+1}.$
\end{theorem}

\dow 
Since $\deg [F,G]<d+i,$ we see by Theorem~\ref{prop_formulaB} that
\begin{equation}\label{row_propHDzieliF3_t}
G_i=\sum_{\substack{\alpha=(\alpha_0,\ldots,\alpha_s)\in\mathbb{Z}\times\mathbb{N}^s\\
t\alpha_0+\alpha_1+2\alpha_2+\cdots+s\alpha_s=i\\
t\alpha_0+d(\alpha_1+\cdots+\alpha_s)\leq N}}
c^{(i)}_{\alpha}F_1^{\alpha_1}\cdots F_s^{\alpha_s}h^{\alpha_0}.
\end{equation}
By Lemma~\ref{lem_MinimalValueAlpha0Plus_kAlpha_s_t} used for $k=0,$ all summands of (\ref{row_propHDzieliF3_t}), except $c^{(i)}_{(\frac{1}{t}(di-sN),0,\ldots,0,N-i)}F_s^{N-i}h^{\frac{1}{t}(di-sN)},$ involves $h^u$ with $u>\tfrac{1}{t}(di-sN).$
Let us also notice that $\tfrac{1}{t}(di-sN)<0.$
Thus, multiplying both sides of (\ref{row_propHDzieliF3_t}) by $h^{\frac{1}{t}(sN-di)},$ we obtain that $h$ divides $c^{(i)}_{(\frac{1}{t}(di-sN),0,\ldots,0,N-i)}F_s^{N-i}.$ 
Since $N\not\equiv 0 (\mod d),$ we see that $-sN\not\equiv 0 (\mod d)$ and so $c^{(i)}_{(\frac{1}{t}(di-sN),0,\ldots,0,N-i)}=\frac{\prod_{k=1}^{N-i}(\tfrac{1}{t}(di-sN)+rk)}{\prod_{i_s=1}^{N-i}ri_s}a_N
=\frac{\prod_{k=1}^{N-i}(di-sN+dk)}{\prod_{i_s=1}^{N-i}di_s}a_N
\neq 0.$
Thus, we see that $h|F_s^{N-i}$ and so $h|F_s$ because $h$ is a square-free polynomial.

Now, we can (\ref{row_propHDzieliF3_t}) rewrite in the form
\begin{equation}\label{row_propHDzieliF3_t2}
G_i=\sum_{\substack{\alpha=(\alpha_0,\ldots,\alpha_s)\in\mathbb{Z}\times\mathbb{N}^s\\
t\alpha_0+\alpha_1+2\alpha_2+\cdots+s\alpha_s=i\\
t\alpha_0+d(\alpha_1+\cdots+\alpha_s)\leq N}}
c^{(i)}_{\alpha}F_1^{\alpha_1}\cdots F_{s-1}^{\alpha_{s-1}}\left(\tfrac{F_s}{h}\right)^{\alpha_s}h^{\alpha_0+\alpha_s},
\end{equation}
where $\tfrac{F_s}{h}$ is a homogeneous polynomial of degree $s-t.$
Using Lemma~\ref{lem_MinimalValueAlpha0Plus_kAlpha_s_t} for $k=1,$ we see that all summands of (\ref{row_propHDzieliF3_t2}), except $c^{(i)}_{(\frac{1}{t}(di-sN),0,\ldots,0,N-i)}\left(\tfrac{F_s}{h}\right)^{N-i}h^{\frac{1}{t}[(d-1)i-(s-1)N]},$ involves $h^u$ with $u>\tfrac{1}{t}[(d-1)i-(s-1)N].$
Let us also notice that $\frac{1}{t}[(d-1)i-(s-1)N]<0.$
Thus, multiplying both sides of (\ref{row_propHDzieliF3_t2}) by $h^{\frac{1}{t}[(d-1)i-(s-1)N]},$ we obtain that $h$ divides $c^{(i)}_{(\frac{1}{t}(di-sN),0,\ldots,0,N-i)}\left(\tfrac{F_s}{h}\right)^{N-i}.$ 
Since $c^{(i)}_{(\frac{1}{t}(di-sN),0,\ldots,0,N-i)}\neq 0,$ it follows as above that $h$ divides $\tfrac{F_s}{h}.$

Now, we can (\ref{row_propHDzieliF3_t2}) rewrite in the form
\begin{equation}\label{row_propHDzieliF3_t3}
G_i=\sum_{\substack{\alpha=(\alpha_0,\ldots,\alpha_s)\in\mathbb{Z}\times\mathbb{N}^s\\
t\alpha_0+\alpha_1+2\alpha_2+\cdots+s\alpha_s=i\\
t\alpha_0+d(\alpha_1+\cdots+\alpha_s)\leq N}}
c^{(i)}_{\alpha}F_1^{\alpha_1}\cdots F_{s-1}^{\alpha_{s-1}}\left(\tfrac{F_s}{h^2}\right)^{\alpha_s}h^{\alpha_0+2\alpha_s},
\end{equation}
where $\tfrac{F_s}{h^2}$ is a homogeneous polynomial of degree $s-2t.$
Thus, we can repeat the above arguments once more again.

Repeating this arguments several times we obtain the thesis.
$\rightBox$ 

\begin{corollary}
	Let $F$ and $G$ be as in Theorem \ref{prop_formulaB}.
	Moreover, let $h$ be a square-free polynomial.
	If $N\not\equiv 0\, (\mod d)$ and $\deg[F,G]< d+i$ with $i<\tfrac{s-k}{d-k}N$ and $2kt<d,$ then we have 
	\begin{equation}
	G_j=\sum_{\substack{\alpha=(\alpha_0,\ldots,\alpha_s)\in\mathbb{Z}\times\mathbb{N}^s\\
			t\alpha_0+\alpha_1+2\alpha_2+\cdots+s\alpha_s=j\\
			t\alpha_0+d(\alpha_1+\cdots+\alpha_s)\leq N}}
	c^{(j)}_{\alpha}F_1^{\alpha_1}\cdots F_{s-1}^{\alpha_{s-1}}\tilde{F}_{s-t(k+1)}^{\alpha_s}h^{\alpha_0+t(k+1)\alpha_s},
	\end{equation}
	where $a_i,\ldots,a_{N-1}\in \mathbb{C}$ and $c^{(j)}_{\alpha}$ are as in Theorem \ref{prop_formulaB} and $\tilde{F}_{s-t(k+1)}$ is as in Theorem~\ref{prop_divisibility_Fs_h2_t}.
\end{corollary}

\section{Dependence between $F_{s-1}$ and $F_s$}\label{sect_dependence_Fs_Fs_minus_1}
In this section we consider the situation when $d_1=\gcd(d,N)=\gcd(\deg F,\deg G)$ is such that $d=2d_1.$ In such a situation $N$ must be an odd multiplicity of $d_1$ and so must be of the form $N=d_1(2k+1)$ for some nonnegative integer $k.$
Since we always assume that $N\geq d,$ we see that the number $k$ maybe considered as positive integer.
In this setup we will show that if $\deg [F,G]$ is small enough that there is a relation between $F_{s-1}$ and $F_s.$ 

We start this with the following numerical lemma which we will use in both considered cases, namely the case $\deg h=1$ and the case $\deg h>1.$

\begin{lemma}\label{lem_Z_Sets_and_Delta_Simplex}
For a given positive integers $j,N,d,t$ such that $j<N,$ $d\leq N$ and $t|\gcd(d,N)$ let us consider the following sets
\begin{eqnarray}
Z_{j,N,d}
&=&
\left\{ (\alpha_0,\ldots,\alpha_s)\in\mathbb{Z}\times\mathbb{N}^s \,|\, \alpha_0+\sum_{l=1}^sl\alpha_l=j, \alpha_0+d\sum_{l=1}^s\alpha_l\leq N \right\},
\\
Z_{j,N,d}^{(t)}
&=&
\left\{ (\alpha_0,\ldots,\alpha_s)\in\mathbb{Z}\times\mathbb{N}^s \,|\, t\alpha_0+\sum_{l=1}^sl\alpha_l=j, t\alpha_0+d\sum_{l=1}^s\alpha_l\leq N \right\},
\\
\Delta_{j,N,d}
&=&
\left\{ (\alpha_1,\ldots,\alpha_s)\in\mathbb{R}_+^s \,|\, \sum_{l=1}^s(d-l)\alpha_l\leq N-j \right\}.
\end{eqnarray}
Then, the following assertions hold:\vspace{0.5cm}\newline
(a)\hspace{0.25cm}
\begin{minipage}{14cm}
The mapping $\varphi:\Delta_{j,N,d}\cap\mathbb{N}^s\ni(\alpha_1,\ldots,\alpha_s)\mapsto (j-\sum_{l=1}^s l\alpha_l,\alpha_1,\ldots,\alpha_s)\in Z_{j,N,d}$ is a bijection.
\end{minipage}
\newline
(b)\hspace{0.25cm}
\begin{minipage}{14cm}
The mapping $\psi:Z_{j,N,d}\cap(t\mathbb{Z}\times\mathbb{N}^s)\ni(\tilde{\alpha}_0,\alpha_1,\ldots,\alpha_s)\mapsto (\frac{\tilde{\alpha}_0}{t},\alpha_1,\ldots,\alpha_s)\in Z_{j,N,d}^{(t)}$ is a bijection.
\end{minipage}
\newline
(c)\hspace{0.25cm}
\begin{minipage}{14cm}
The set $\Delta_{j,N,d}$ is the simplex with the following set of vertices: $(0,\ldots,0),$ $\frac{N-j}{s}e_1,$ $\frac{N-j}{s-1}e_2,\ldots,$ $\frac{N-j}{2}e_{s-1}$ and $(N-j)e_s,$ where $e_1=(1,0,\ldots,0),\ldots,e_s=(0,\ldots,0,1)$ is the standard basis of $\mathbb{R}^s.$
In other words vertices of $\Delta_{j,N,d}$ others than $(0,\ldots,0)$ are of the form $\frac{N-j}{d-l}e_l$ for $l=1,\ldots,s.$
\end{minipage}
\newline
(d)\hspace{0.25cm}
\begin{minipage}{14cm}
For any $(\beta_1,\ldots,\beta_s)\in\mathbb{N}^s$ and for the function (a linear function actually)
$
f_{\beta_1,\ldots,\beta_s}:\Delta_{j,N,d}\ni(\alpha_1,\ldots,\alpha_s)\mapsto j-\sum_{l=1}^s
(l-\beta_l)\alpha_l\in\mathbb{R}
$
we have $f_{\beta_1,\ldots,\beta_s}(0,\ldots,0)=j$ and $f_{\beta_1,\ldots,\beta_s}(\frac{N-j}{d-l}e_l)=j-\frac{l-\beta_l}{d-l}(N-j).$
\end{minipage}
\newline
(e)\hspace{0.25cm}
\begin{minipage}{14cm}
If $2\beta_s<d$ then the minimal value of $f_{0,\ldots,0,\beta_s}$ on the set $\Delta_{j,N,d}$ equals to $j-(s-\beta_s)(N-j)$ and is attained exactly at one point $(0,\ldots,0,N-j).$
\end{minipage}
\newline
(f)\hspace{0.25cm}
\begin{minipage}{14cm}
If $2\beta_s=d$ then the minimal value of $f_{0,\ldots,0,\beta_s}$ on the set $\Delta_{j,N,d}$ equals to $j-(s-\beta_s)(N-j)$ and is attained exactly on the line segment with the following endpoints $(0,\ldots,0,\frac{N-j}{2},0)$ and $(0,\ldots,0,N-j).$
\end{minipage}
\end{lemma}
\dow
The assertions (a)-(d) are easy to check.
To prove assertions (e) and (f) let us notice that, since the set $\Delta_{j,N,d}$ is a convex set and the function $f_{0,\ldots,0,\beta_s}$ is linear, it follows that we only need to examine the values of the function on the set of vertices of $\Delta_{j,N,d}.$ 
Since we are interested in minimal value of $f_{0,\ldots,0,\beta_s}$ and $N-j>0,$ we see that what is interested for us is the maximal value of the fraction $\frac{l-\beta_l}{d-l}$ for $l=1,\ldots,s.$
Since we have $\beta_1=\ldots=\beta_{s-1}=0,$ the values of $\frac{l-\beta_l}{d-l}$ for $l=1,\ldots,s-1$ are the following: $\frac{1}{d-1}<\frac{2}{d-2}<\ldots<\frac{d-2}{2}=\frac{s-1}{2}.$
To conclude the proof, we only need to notice that $\frac{s-1}{2}< s-\beta_s$ in the case $2\beta_s<d$ and $\frac{s-1}{2}=s-\beta_s$ in the case $2\beta_s=d.$
$\rightBox$

Before we go to the subsection devoted to the case $\deg h=1$ we give the one more lemma, which is a consequence of the above numerical one.

\begin{lemma}\label{lem_MinimalValuesSets}
Let $d_1\in\mathbb{N}\setminus\{0,1\}$ and set $d=2d_1,$ $N=d_1(2k+1)$ for some $k\in\mathbb{N}_+.$
Then, the following statements hold \vspace{0.2cm}\newline
(a)\hspace{0.25cm}
\begin{minipage}{14cm}
$\min \{\, \alpha_0+d_1\alpha_s\,|\, (\alpha_0,\ldots,\alpha_s)\in Z_{j,N,d}\, \}$ is equal to $(d-d_1)j-(s-d_1)N=d_1j-(d_1-1)N$ and is attained exactly at the points of the following set 
$\{ \,(dj-sN+dl,0,\ldots,0,l,N-j-2l)\in\mathbb{Z}\times\mathbb{N}^s \, 
|\, 0\leq 2l\leq N-j\, \},$
\end{minipage}
\vspace{0.2cm}\newline
(b)\hspace{0.25cm}
\begin{minipage}{14cm}
Let $t$ be an integer such that $t$ divides $d_1$ and set $r=\frac{d}{t}$ and $\tilde{r}=\frac{d_1}{t}$ (in particular $r=2\tilde{r}$).
Then, $\min \{\, \alpha_0+\tilde{r}\alpha_s\,|\, (\alpha_0,\ldots,\alpha_s)\in Z_{j,N,d}^{(t)}\, \}$ is equal to $\tilde{r}j-(d_1-1)\tilde{r}(2k+1)$ and is attained exactly at the points of the following set 
$\{ \,(rj-s\tilde{r}(2k+1)+rl,0,\ldots,0,l,N-j-2l)\in\mathbb{Z}\times\mathbb{N}^s \, 
|\, 0\leq 2l\leq N-j\, \}.$
\end{minipage}
\vspace{0.2cm}\newline
\end{lemma}
\dow
Consider the set $\Delta_{j,N,d}
=\left\{ (\alpha_1,\ldots,\alpha_s)\in\mathbb{R}_+^s \,|\, \sum_{l=1}^s(d-l)\alpha_l\leq N-j \right\}$ and the function
$ f_{0,\ldots,d_1}(\alpha_1,\ldots,\alpha_s) = 
j-\sum_{l=1}^s
(l-\delta_{ls}d_1)\alpha_l,$ where $\delta_{ls}$ is the Kronecker delta.
Let us notice that from Lemma~\ref{lem_Z_Sets_and_Delta_Simplex} (a), we have $Z_{j,N,d}=\varphi_j(\Delta_{j,N,d}\cap\mathbb{N}^s).$
Notice also that  $f_{0,\ldots,d_1}(\alpha_1,\ldots,\alpha_s)=\alpha_0+d_1\alpha_s,$ 
where $(\alpha_0,\alpha_1,\ldots,\alpha_s)=\varphi_j(\alpha_1,\ldots,\alpha_s)$ 
for any $(\alpha_1,\ldots,\alpha_s) \in \Delta_{j,N,d}\cap\mathbb{N}^s.$
Now using the above considerations and Lemma~\ref{lem_Z_Sets_and_Delta_Simplex} (c) and (f), we obtain the assertion (a).

In order to prove (b) we apply Lemma~\ref{lem_Z_Sets_and_Delta_Simplex} (b) and already proved Lemma \ref{lem_MinimalValuesSets} (a).
Indeed, by Lemma~\ref{lem_Z_Sets_and_Delta_Simplex} (b) we have that $(\alpha_0,\ldots,\alpha_s)\in Z_{j,N,d}^{(t)}$ if and only if  $(\tilde{\alpha}_0,\ldots,\alpha_s)=(t\alpha_0,\ldots,\alpha_s) \in Z_{j,N,d}\cap(t\mathbb{Z}\times\mathbb{N}^s).$
On the other hand minimizing $\alpha_0+\tilde{r}\alpha_s$ is equivalent to minimizing $t(\alpha_0+\tilde{r}\alpha_s)=\tilde{\alpha}_0+d_1\alpha_s.$
Thus searching the points at which $\alpha_0+\tilde{r}\alpha_s$ achieve minimal values on the set $Z_{j,N,d}^{(t)}$ is the same as searching the points at which $\tilde{\alpha}_0+d_1\alpha_s$ achieve minimal values on the set $Z_{j,N,d}\cap(t\mathbb{Z}\times\mathbb{N}^s)$ and so we can use Lemma \ref{lem_MinimalValuesSets} (a). 
It should be mentioned that from Lemma \ref{lem_MinimalValuesSets} (a) we get information about the points at which the minimum is achieved and what is its value.
This completes the proof.
$\rightBox$\\

To conclude this part of the section, we make the following observation.
\begin{remark}
	One can notice that Lemma~\ref{lem_Z_Sets_and_Delta_Simplex} (e) can be used, instead of  Lemmas~\ref{lem_MinimalValueAlpha0} and \ref{lem_MinimalValueAlpha0Plus_kAlpha_s_t}, in order to obtain alternative proofs of Theorems~\ref{prop_F3equal_f2h} and \ref{prop_divisibility_Fs_h2_t}.
\end{remark}

\subsection{The case $\deg h=1$}

In this subsection we consider the case $\deg h=1.$
Thus, we have $F=F_1+\cdots+F_{s-1}+h^d$ and $G=G_1+\cdots+G_{N-1}+a_Nh^N$ with $d=2d_1$ for some $d_1\in\mathbb{N}\setminus\{0,1\},$ $a_N\in\mathbb{C}\setminus\{0\}$ and $N=d_1(2k+1)$ for some positive integer $k.$

In the proof of the main result of this subsection we will use the following lemma.

\begin{lemma}\label{lem_Fs_minus1FsSymbolNewtona}
Assume that $d=2d_1$ for some $d_1\in\mathbb{N}\setminus\{0,1\},$ $N=2d_1k+d_1$ for some $k\in\mathbb{N}_+ ,$ and that $\deg [F,G]<d+N-2k-2.$ 
Assume also that $F_d=h^d$ and $G_N=a_Nh^N$ for some homogeneous polynomial $h$ of degree one and $a_N\in\mathbb{C}\setminus\{0\}.$
Then, for all $l\in\{ 0,1,\ldots, k+1 \} ,$ we have
\begin{equation}
c^{(N-2k-2)}_{(-2d_1k-3d_1+2d_1l,0,\ldots,0,l,2k+2-2l)}=(-4)^l {k+1 \choose l} c^{(N-2k-2)}_{(-2d_1k-3d_1,0,\ldots,0,0,2k+2)}.
\end{equation}
\end{lemma}
\dow 
Since $\deg [F,G]<d+N-2k-2,$ we see that the polynomial $G_{N-2k-2}$ is given by the formulas (\ref{row_formula1}) and (\ref{row_formula2}) with $d=2d_1$ and $s=2d_1-1.$
One can check, that
\begin{eqnarray}
&&
c^{(N-2k-2)}_{(-2d_1k-3d_1+2d_1(l+1),0,\ldots,0,l+1,2k+2-2(l+1))}
\\
\nonumber
&=&
\frac{\prod_{i=1}^{2k+2-(l+1)}(-2d_1k-3d_1+2d_1(l+1)+2d_1i)}{\prod_{i_{s-1}=1}^{l+1}2d_1i_{s-1}\cdot \prod_{i_s=1}^{2k+2-2(l+1)}2d_1i_s}a_N
\end{eqnarray}
and
\begin{eqnarray}
c^{(N-2k-2)}_{(-2d_1k-3d_1+2d_1l,0,\ldots,0,l,2k+2-2l)}
&=&
\frac{\prod_{i=1}^{2k+2-l}(-2d_1k-3d_1+2d_1l+2d_1i)}{\prod_{i_{s-1}=1}^{l}2d_1i_{s-1}\cdot \prod_{i_s=1}^{2k+2-2l}2d_1i_s}a_N,
\end{eqnarray}
for any $l=0,1,\ldots,k.$
From the above, one can see that
\begin{eqnarray}
&&
c^{(N-2k-2)}_{(-2d_1k-3d_1+2d_1(l+1),0,\ldots,0,l+1,2k+2-2(l+1))}
\\
&=&
\nonumber
\frac{2d_1(2k+1-2l)2d_1(2k+2-2l)}{2d_1(l+1)(-2d_1k-d_1+2d_1l)}c^{(N-2k-2)}_{(-2d_1k-3d_1+2d_1l,0,\ldots,0,l,2k+2-2l)}
\\
&=&
\frac{(k+1-l)(-4)}{l+1}c^{(N-2k-2)}_{(-2d_1k-3d_1+2d_1l,0,\ldots,0,l,2k+2-2l)}.\nonumber
\end{eqnarray}
In particular, for $l=0,$ we obtain that
\begin{eqnarray}
c^{(N-2k-2)}_{(-2d_1k-d_1,0,\ldots,0,1,2k)}
&=&
\frac{(k+1)(-4)}{1}c^{(N-2k-2)}_{(-2d_1k-3d_1,0,\ldots,0,0,2k+2)}
\\
&=&
\nonumber
(-4){k+1 \choose 1}c^{(N-2k-2)}_{(-2d_1k-3d_1,0,\ldots,0,0,2k+2)}.
\end{eqnarray}
Similarly, one can see that
\begin{eqnarray}
c^{(N-2k-2)}_{(-2d_1k+d_1,0,\ldots,0,2,2k-2)}
&=&
\frac{(k+1-1)(-4)}{1+1}c^{(N-2k-2)}_{(-2d_1k-d_1,0,\ldots,0,1,2k)}
\nonumber\\
&=&
\frac{(k+1)k(-4)^2}{1\cdot 2}c^{(N-2k-2)}_{(-2d_1k-3d_1,0,\ldots,0,0,2k+2)}
\nonumber
\\
&=&
(-4)^2{k+1 \choose 2}c^{(N-2k-2)}_{(-2d_1k-3d_1,0,\ldots,0,0,2k+2)}
.\nonumber
\end{eqnarray}
Continuing the above calculation inductively, we obtain the result.
$\rightBox$

Now we can prove the following theorem.

\begin{theorem}\label{prop_DependeceF2F3}
Assume that $d=2d_1$ for some $d_1\in\mathbb{N}\setminus\{0,1\},$ $N=2d_1k+d_1$ for some $k\in\mathbb{N}_+ ,$ and that $\deg [F,G]<d+N-2k-2.$ 
Assume also that $F_d=h^d$ and $G_N=a_Nh^N$ for some homogeneous polynomial $h$ of degree one and $a_N\in\mathbb{C}\setminus\{0\}.$ 
Then, there is a homogeneous polynomial $\hat{F}_{s-2}$ of degree $s-2$ such that
\begin{equation}
F_{s-1}=\tfrac{1}{4}(\tilde{F}_{d_1-1}^2+h\hat{F}_{s-2}).
\end{equation}
\end{theorem}
\dow 
First of all notice that $\frac{s-(d_1-1)}{d-(d_1-1)}=\frac{d_1}{d_1+1}> \frac{d_1-1}{d_1}$ and so $\frac{s-(d_1-1)}{d-(d_1-1)}N=\frac{d_1}{d_1+1}N>\frac{d_1-1}{d_1}N=\frac{d_1-1}{d_1}(2d_1k+d_1)=N-2k-1.$
By the above and the assumption that $\deg [F,G]<d+N-2k-2,$ we can use Theorem~\ref{prop_F3equal_f2h2} to obtain that $F_s$ is divisible by $h^{d_1}.$
Thus, there exists a homogeneous polynomial $\tilde{F}_{d_1-1}$ of degree $d_1-1$ such that $F_s=h^{d_1}\tilde{F}_{d_1-1}.$
Now, using the above considerations, assumption that $\deg [F,G]<d+N-2k-2$ and Theorem~\ref{prop_formula} with $d=2d_1,$ $s=2d_1-1$ and $j=N-2k-2$ we obtian that
\begin{equation}\label{row_PropDependenceF2F3}
G_{N-2k-2}=\sum_{\substack{\alpha=(\alpha_0,\ldots,\alpha_s)\in\mathbb{Z}\times\mathbb{N}^s\\
\alpha_0+\alpha_1+2\alpha_2+\cdots+s\alpha_s=j\\
\alpha_0+d(\alpha_1+\cdots+\alpha_s)\leq N}}
c^{(N-2k-2)}_{\alpha}F_1^{\alpha_1}\cdots F_{s-1}^{\alpha_{s-1}}\tilde{F}_{d_1-1}^{\alpha_s}h^{\alpha_0+d_1\alpha_s}.
\end{equation}

Notice that, for $N=2d_1k+d_1$ and $j=N-2k-2$, the set from Lemma~\ref{lem_MinimalValuesSets}~(a) is equal to the following one 
$\{ \,(-2d_1k-3d_1+2d_1l,0,\ldots,0,l,2k+2-2l)\in\mathbb{Z}\times\mathbb{N}^s \, 
|\, l=0,1,\ldots,k+1 \, \}.$
Thus, in the polynomial $G_{N-2k-2}$
all summand other than included in the following sum
\begin{equation}
\sum_{l=0}^{k+1}c^{(N-2k-2)}_{(-2d_1k-3d_1+2d_1l,0,\ldots,0,l,2k+2-2l)}F_{s-1}^{l}\tilde{F}_{d_1-1}^{2k+2-2l}h^{-d_1}
\end{equation}
involves $h^u$ with $u>-d_1.$ 
By Lemma~\ref{lem_Fs_minus1FsSymbolNewtona}, the above sum can be written as follows
\begin{eqnarray}
&&
c^{(N-2k-2)}_{(-2d_1k-3d_1,0,\ldots,0,2k+2)}h^{-d_1}\sum_{l=0}^{k+1}{k+1 \choose l} (-4)^l F_{s-1}^l\tilde{F}_{d_1-1}^{2k+2-2l}
\\
\nonumber
&=& 
c^{(N-2k-2)}_{(-2d_1k-3d_1,0,\ldots,0,2k+2)}h^{-d_1} (\tilde{F}_{d_1-1}^2-4F_{s-1})^{k+1}.
\end{eqnarray}

Now, we see by multiplying both sides of (\ref{row_PropDependenceF2F3}) by $h^{d_1}$ that $h$ divides $c^{(N-2k-2)}_{(-2d_1k-3d_1,0,\ldots,0,2k+2)} (\tilde{F}_{d_1-1}^2 \allowbreak -4F_{s-1})^{k+1}.$ 
But, since $c^{(N-2k-2)}_{(-2d_1k-3d_1,0,\ldots,0,2k+2)}=\frac{\prod_{i=1}^{2k+2}(-2d_1k-3d_1+2d_1i)}{\prod_{i_s=1}^{2k+2}2d_1i_s}a_{2d_1k+d_1}\neq 0,$ it follows that $h|(\tilde{F}_{d_1-1}^2-4F_{s-1})^{k+1},$ and since $h$ is not a power of any other polynomial (in this case this means that $h$ is square-free), $h|(\tilde{F}_{d_1-1}^2-4F_{s-1}).$
Thus there is a homogeneous polynomial $\hat{F}_{s-2}$ of degree $s-2$ such that $-h\hat{F}_{s-2}=\tilde{F}_{d_1-1}^2-4F_{s-1}.$ Solving the last equality with respect to $F_{s-1}$ we obtain the thesis. 
$\rightBox$

\subsection{The case $\deg h>1$}
As the title of subsection suggests, we consider here the situation of $\deg h>1.$
Thus, we have $F=F_1+\cdots+F_{s-1}+h^{\frac{d}{t}}$ and $G=G_1+\cdots+G_{N-1}+a_Nh^{\frac{N}{t}}$ where $t=\deg h$ and $a_N\in\mathbb{C}\setminus\{0\}.$
As before, we assume that $d=2d_1$ for some $d_1\in\mathbb{N}\setminus\{0,1\}$  and $N=d_1(2k+1)$ for some positive integer $k.$ 
Since $t$ must divide $d$ and $N$ and $\gcd(d,N)=d_1,$ one can see that we must have, in this situation, that $d_1$ is divisible by $t.$

In the proof of the main result of this subsection we will use the following lemma.

\begin{lemma}\label{lem_Fs_minus1FsSymbolNewtona_t}
	Let $d, N$ and $\deg [F,G]$ be as in Lemma \ref{lem_Fs_minus1FsSymbolNewtona}.
	Moreover, suppose that $F_d=h^\frac{d}{t}$ and $G_N=a_Nh^\frac{N}{t}$ for some square-free, homogeneous polynomial $h$ of degree $t$, where $t \geq 1,$ $a_N\in\mathbb{C}\setminus\{0\}$ and $r=\frac{d}{t}$, $\tilde{r}=\frac{d_1}{t}$.
	Then, for all $l\in\{ 0,1,\ldots, k+1 \} ,$ we have
	\begin{equation}
	c^{(N-2k-2)}_{((-2k-3+2l)\tilde{r},0,\ldots,0,l,2k+2-2l)}=(-4)^l {k+1 \choose l} c^{(N-2k-2)}_{((-2k-3)\tilde{r},0,\ldots,0,0,2k+2)}.
	\end{equation}
\end{lemma}
\dow 
According with the assumption that $\deg [F,G]<d+N-2k-2,$ the polynomial $G_{N-2k-2}$ is given by the formulas (\ref{row_formula1B}) and (\ref{row_formula2B}) with $d=2d_1$, $s=2d_1-1,$ so
\begin{eqnarray}
&&
c^{(N-2k-2)}_{((-2k-3+2(l+1))\tilde{r},0,\ldots,0,l+1,2k+2-2(l+1))}
\\
\nonumber
&=&
c^{(N-2k-2)}_{\left(\frac{-2d_1k-3d_1+2d_1(l+1)}{t},0,\ldots,0,l+1,2k+2-2(l+1)\right)}
\\
\nonumber
&=&
\frac{\prod_{i=1}^{2k+2-(l+1)}\frac{-2d_1k-3d_1+2d_1(l+1)+2d_1i}{t}}{\prod_{i_{s-1}=1}^{l+1}\frac{2d_1i_{s-1}}{t}\cdot \prod_{i_s=1}^{2k+2-2(l+1)}\frac{2d_1i_s}{t}}a_N
\\
\nonumber
&=&
\frac{\prod_{i=1}^{2k+2-(l+1)}(-2d_1k-3d_1+2d_1(l+1)+2d_1i)}{\prod_{i_{s-1}=1}^{l+1}2d_1i_{s-1}\cdot \prod_{i_s=1}^{2k+2-2(l+1)}2d_1i_s}a_N
\end{eqnarray}
and
\begin{eqnarray}
c^{(N-2k-2)}_{((-2k-3+2l)\tilde{r},0,\ldots,0,l,2k+2-2l)}
&=&
\frac{\prod_{i=1}^{2k+2-l}(-2d_1k-3d_1+2d_1l+2d_1i)}{\prod_{i_{s-1}=1}^{l}2d_1i_{s-1}\cdot \prod_{i_s=1}^{2k+2-2l}2d_1i_s}a_N,
\end{eqnarray}
for any $l=0,1,\ldots,k.$ 
Further part of this proof is similar to the proof of Lemma \ref{lem_Fs_minus1FsSymbolNewtona}.
$\rightBox$

\begin{theorem}\label{prop_DependeceF2F3, degh>1}
	Assume that $d=2d_1$ for some $d_1\in\mathbb{N}\setminus\{0,1\},$ $N=2d_1k+d_1$ for some $k\in\mathbb{N}_+ ,$ and that $\deg [F,G]<d+N-2k-2.$ 
	Assume also that $F_d=h^\frac{d}{t}=h^r$ and $G_N=a_Nh^\frac{N}{t}$ for some square-free, homogeneous polynomial $h$ of degree $t$, where $t\geq1$, $a_N\in\mathbb{C}\setminus\{0\}$ and $r=\frac{d}{t}$, $\tilde{r}=\frac{d_1}{t}$.
	Then, there is a homogeneous polynomial $\hat{F}_{s-t-1}$ of degree $s-t-1$ such that
	\begin{equation}
	F_{s-1}=\tfrac{1}{4}(\tilde{F}_{d_1-1}^2+h\hat{F}_{s-t-1}).
	\end{equation}
\end{theorem}
\dow 
Using the assumption that $\deg [F,G]<d+N-2k-2,$ and the following observation that
$\frac{s-(\frac{d_1}{t}-1)}{d-(\frac{d_1}{t}-1)}=\frac{d_1(2t-1)}{d_1(2t-1)+t}> \frac{d_1-1}{d_1}$ and so $\frac{s-(\frac{d_1}{t}-1)}{d-(\frac{d_1}{t}-1)}N>\frac{d_1-1}{d_1}N=N-2k-1$
we can apply Theorem \ref{prop_divisibility_Fs_h2_t} to obtain the divisibility of $F_s$ by $h^{\tilde{r}}.$
It means that, there exists a homogeneous polynomial $\tilde{F}_{d_1-1}$ of degree $d_1-1$ such that $F_s=h^{\tilde{r}}\tilde{F}_{d_1-1}.$
Referring to the above calculations and Theorem \ref{prop_formulaB}, the polynomial $G_{N-2k-2}$ takes the form
\begin{equation}\label{row_PropDependenceF2F31}
G_{N-2k-2}=\sum_{\substack{\alpha=(\alpha_0,\ldots,\alpha_s)\in\mathbb{Z}\times\mathbb{N}^s\\
		t\alpha_0+\alpha_1+2\alpha_2+\cdots+s\alpha_s=j\\
		t\alpha_0+d(\alpha_1+\cdots+\alpha_s)\leq N}}
c^{(N-2k-2)}_{\alpha}F_1^{\alpha_1}\cdots F_{s-1}^{\alpha_{s-1}}\tilde{F}_{d_1-1}^{\alpha_s}h^{\alpha_0+\tilde{r}\alpha_s}.
\end{equation}

According to Lemma \ref{lem_MinimalValuesSets}~(b)
for $j=N-2k-2$, the set of summands with minimal possible power of $h$, in the above sum, corresponds to $\alpha$ belonging to
$\{ \,((-2k-3+2l)\tilde{r},0,\ldots,0,l,2k+2-2l)\in\mathbb{Z}\times\mathbb{N}^s \, 
|\, l=0,1,\ldots,k+1 \, \},$
so
all components under than occurring the presenting sum
\begin{equation}
\sum_{l=0}^{k+1}c^{(N-2k-2)}_{((-2k-3+2l)\tilde{r},0,\ldots,0,l,2k+2-2l)}F_{s-1}^{l}\tilde{F}_{d_1-1}^{2k+2-2l}h^{-\tilde{r}}
\end{equation}
involves $h^u$ with $u>-\tilde{r}.$ 
Applying Lemma~\ref{lem_Fs_minus1FsSymbolNewtona_t}, we can rewrite the above sum in the following way
\begin{eqnarray}
&&
c^{(N-2k-2)}_{((-2k-3)\tilde{r},0,\ldots,0,2k+2)}h^{-\tilde{r}}\sum_{l=0}^{k+1}{k+1 \choose l} (-4)^l F_{s-1}^l\tilde{F}_{d_1-1}^{2k+2-2l}
\\
\nonumber
&=& 
c^{(N-2k-2)}_{((-2k-3)\tilde{r},0,\ldots,0,2k+2)}h^{-\tilde{r}} (\tilde{F}_{d_1-1}^2-4F_{s-1})^{k+1}.
\end{eqnarray}

The rest of the proof is analogous to the proof of the Theorem \ref{prop_DependeceF2F3} using the assumption that $h$ is square-free.
$\rightBox$




\vspace{1cm}
 
\textsc{Daria Holik\newline
AGH University of Science and Technology,\newline
Faculty of Applied Mathematics\newline
al. A. Mickiewicza 30}\newline
\textsc{30-059 Krak\'{o}w\newline
Poland\newline
} e-mail: holikd@agh.edu.pl \newline
ORCID: https://orcid.org/0000-0002-2133-9711
\vspace{0.5cm}

\textsc{Marek Kara\'{s}\newline
AGH University of Science and Technology,\newline
Faculty of Applied Mathematics\newline
al. A. Mickiewicza 30}\newline
\textsc{30-059 Krak\'{o}w\newline
Poland\newline
} e-mail: mkaras@agh.edu.pl\newline
ORCID: https://orcid.org/0000-0003-0821-521X
\vspace{0.5cm}


\begin{thebibliography}{99}


\bibitem{Drensky_Yu}
 V. Drensky, J.-T. Yu, 
 \tit{Degree estimate for commutators.} 
 J. Algebra \notextbf{322}, no. 7, 2321–2334 (2009)

\bibitem{EdoKanKarKur}
 E. Edo, T. Kanehira, M. Karaś,  S. Kuroda, 
 \tit{Separability of wild automorphisms of a polynomial ring.} Transform. Groups \notextbf{18}, no. 1, 81–96 (2013)

\bibitem{Jedrzejewicz_Zielinski}
 P. J\k{e}drzejewicz, J. Zieliński,  
 \tit{Analogs of Jacobian conditions for subrings.} 
 J.\ Pure\ Appl.\ Algebra \notextbf{221}, 2111-2118 (2017)

\bibitem{Karas1} 
 M. Karaś, 
 \tit{There is no tame automorphism of} $ \mathbb{C}^{3}$ \tit{with multidegree} $(3,4,5).$ 
 Proc. Am. Math. Soc. \notextbf{139}, no. 3, 769-775 (2011)

\bibitem{Karas3} 
 M. Karaś,
 \tit{Tame automorphisms of }$\mathbb{C}^{3} $ \tit{with multidegree of the form }$(3,d_{2},d_{3}).$ 
 J.\ Pure\ Appl.\ Algebra \notextbf{214}, no. 12, 2144-2147 (2010)

\bibitem{Karas4}
 M. Karaś,
 \tit{Multidegrees of tame automorphisms of} $\mathbb{C}^n.$
 Diss. Math. \notextbf{477}, 55 p. (2011)

\bibitem{KarasZygad}
 M. Karaś, J. Zygadło, 
 \tit{On multidegree of tame and wild automorphisms of} $\mathbb{C}^{3}.$ 
 J.\ Pure\ Appl.\ Algebra \notextbf{215}, 2843–2846 (2011)

\bibitem{Kuroda3}
 S. Kuroda, 
 \tit{A generalization of the Shestakov-Umirbaev inequality.}
 J. Math. Soc. Japan \notextbf{60}, no. 2, 495–510 (2008)

\bibitem{Kuroda2}
 S. Kuroda,
 \tit{Shestakov-Umirbaev reductions and Nagata's conjecture on a polynomial automorphisms.}
 Tohoku Math. J. \notextbf{62}, 75-115 (2010)

\bibitem{Kuroda4} 
 S. Kuroda, 
 \tit{On the Karaś type theorems for the multidegrees of polynomial automorphisms.}
 J. Algebra \notextbf{423}, 441–465 (2015)

\bibitem{Makar_Yu} 
 L. Makar-Limanov, J.-T. Yu,
 \tit{Degree estimate for subalgebras generated by two elements.} 
 J. Eur. Math. Soc. (JEMS) \notextbf{10}, 533–541 (2008) 

\bibitem{Pritchard} 
 F. L. Pritchard, 
 \tit{Polynomial mappings with Jacobian determinant of bounded degree.} 
 Arch. Math. \notextbf{48}, no. 6, 495–504 (1987)

\bibitem{shUmb1} 
 I. P. Shestakov, U. U. Umirbaev,
 \tit{The Nagata automorphism is wild.} 
 Proc. Natl. Acad. Sci. USA \notextbf{100}, 12561-12563 (2003)
 
 \bibitem{shUmb} 
 I. P. Shestakov, U. U. Umirbaev,
 \tit{Poisson brackets and two-generated subalgebras of rings of polynomials.} 
 J. Amer. Math. Soc. \notextbf{17}, 181-196 (2004)

\bibitem{shUmb2}  
 I. P. Shestakov, U. U. Umirbaev,
 \tit{The tame and the wild automorphisms of polynomial rings in three variables.}
 J. Amer. Math. Soc. \notextbf{17}, 197-227 (2004)

\bibitem{SunChen} 
 X. Sun, Y. Chen, 
 \tit{Multidegrees of tame automorphisms in dimension three.} 
 Publ. Res. Inst. Math. Sci. \notextbf{48}, no. 1, 129–137 (2012)


\bibitem{Umirbaev Yu}  
 U. U. Umirbaev, J.-T.Yu,
 \tit{The strong Nagata conjecture.}
 Proc. Natl. Acad. Sci. USA \notextbf{101}, 4352-4355 (2004) 

\bibitem{Yu} 
 J.-T. Yu, 
 \tit{Degree estimate for subalgebras and automorphisms of free algebras.} 
 In: S.-T. Yau, et al. (Eds.), Proceedings of the Fourth International Congress of Chinese Mathematicians, Hangzhou, China, December 17–23, (2007), Higher Education Press, Beijing (2008), pp. 359–366; International Press, Boston (2009)

\end{thebibliography}
\end{document}